\documentclass{article}
\usepackage{graphicx} 
\usepackage{algorithm}
\usepackage{algpseudocode}
\usepackage{amsmath}
\usepackage{amssymb}
\usepackage{tikz}
\usepackage[utf8]{inputenc}
\usepackage{cite}
\usepackage[T1]{fontenc}
\DeclareUnicodeCharacter{0301}{\'{e}}
\DeclareUnicodeCharacter{0308}{\"{o}}
\usepackage{xcolor,hyperref}
\hypersetup{colorlinks,breaklinks,
            linkcolor=blue,urlcolor=blue,
            anchorcolor=blue,citecolor=blue}
\usepackage{todonotes}
\usepackage{array}
\usepackage{booktabs}

\newcommand{\md}[1]{\textcolor{black}{#1}}

\newcommand{\tetrahedron}{%
  \mathrel{\vcenter{\hbox{%
    \begin{tikzpicture}[scale=0.11]
      \draw (0,0) -- (2,0) -- (1,1.732) -- cycle; 
      \draw (1,1.732) -- (1,0.577); 
      \draw (0,0) -- (1,0.577);
      \draw (2,0) -- (1,0.577);
    \end{tikzpicture}%
  }}}%
}

\title{Alternative set-theoretical algorithms for efficient\\ computations of cliques in Vietoris-Rips complexes}
\author{Danillo Barros de Souza\thanks{Basque Center for Applied Mathematics, Bilbao, Spain}, Jonatas Teodomiro\thanks{Universidade Federal de Pernambuco, Recife, Brasil}, Fernando A N Santos\thanks{DIEP, Amsterdam, Netherlands},\\ Mathieu Desroches\thanks{Inria, Montpellier, France} \,and Serafim Rodrigues\thanks{Ikerbasque Science Foundation, Basque Center for Applied Mathematics and ikerbasque, Bilbao, Spain}}
\date{\md{\today}}

\begin{document}

\maketitle
\begin{abstract}
{Identifying cliques in dense networks remains a formidable challenge, even with significant advances in computational power and methodologies. To tackle this, numerous algorithms have been developed to optimize time and memory usage, implemented across diverse programming languages. Yet, the inherent NP-completeness of the problem continues to hinder performance on large-scale networks, often resulting in memory leaks and slow computations. In the present study, we critically evaluate classic algorithms to pinpoint computational bottlenecks and introduce novel set-theoretical approaches tailored for network clique computation. Our proposed algorithms are rigorously implemented and benchmarked against existing Python-based solutions, demonstrating superior performance. These findings underscore the potential of set-theoretical techniques to drive substantial performance gains in network analysis.}

\end{abstract}
\section{Introduction}
Network science \cite{posfai2016network} has emerged as a field providing a powerful framework of tools and methodologies to study complex networks. Its applications and importance span multiple real-world domains, to name a few (but not exclusive), community detection in random graphs~\cite{bollobas1976cliques}, social networks~\cite{fortunato2010community}, protein-protein interactions~\cite{palla2005uncovering}, the detection of co-authorship collaboration networks~\cite{newman2001structure}, and economic indicators~\cite{schweitzer2009economic}. One critical aim in the analysis of complex networks is that of identifying cliques, which is known for its NP-hard complexity~\cite{gary1979computers,boix2021average}. In classical algorithms, a typical strategy in identifying cliques is to build upon lower dimensional cliques, such as triangle cliques formed by a combination of edges ~\cite{chiba1985arboricity,al2018triangle}. More contemporary algorithms also compute higher-order cliques (i.e. cliques mediated by structures beyond pairwise interactions). These can be used to determine maximal cliques \cite{bron1973algorithm,tsukiyama1977new,chiba1985arboricity,tomita2003efficient,eppstein2010listing} and cliques of all dimensions \cite{zhang2005genome,tomita2006worst}. However, the most prominent use of higher-order cliques is via recent developments in topological data analysis (TDA)\cite{pascucci2010topological,zomorodian2012topological}. In this context, \md{clique computation} is performed \md{in the background,} and this includes cliques emerging from pairwise weighted graph structures, the so-called clique complexes. This information is then used to compute topological invariants, such as, Betti numbers~\cite{horak2009persistent,edelsbrunner2014computational}, the \md{barcode} diagrams~\cite{ghrist2008barcodes,bubenik2015statistical}, the Euler characteristics \cite{knill2011dimension,lawniczak2021euler,de2022euler} and the higher-order Forman-Ricci curvature \cite{Barros_de_Souza_2021,de2023efficient}. The computation of high-dimensional cliques is also fundamental in geometrical approaches, for instance, in geometrical data analysis via discrete Ricci curvatures \cite{forman2003bochner,bourne2018ollivier,samal2018comparative}, as well as, when computing Laplacians on networks \cite{lim2020hodge}. These geometrical methods are fundamental across various fields of application, such as neuroscience~\cite{santos2018topological,chatterjee2021detecting,santos2023emergence}, market fragility detection \cite{sandhu2016ricci}, or epidemiology where geometrical data analysis can be used to identify markers and predict epidemic outbreaks~\cite{Barros_de_Souza_2021,de2022euler}. Numerous algorithms have been proposed to address the computational complexity of finding cliques \cite{gary1979computers,chiba1985arboricity,tomita2003efficient,schank2005finding,zhang2005genome,boix2021average,bauer2017phat}, each varying in terms of implementation feasibility and efficiency. However, many of these algorithms face significant challenges, including excessive runtimes for clique detection and memory leaks that hinder their ability to handle large networks. As a result, these algorithms often fall short of efficiently identifying cliques at scale. Set-theoretical approaches \cite{enderton1977elements,jech2003set,zimmermann1985applications,pawlak2002rough} offer valuable alternative theoretical formulations for addressing discrete problems, particularly in the construction of posets \cite{bloch2014combinatorial} and simplicial complexes \cite{yadav2022poset}. For example, our recent work in \cite{de2023efficient} leveraged these methods to provide novel theoretical derivations for efficiently computing the discrete Forman-Ricci curvature in simplicial complexes. The underlying idea was to utilize node adjacency to map and classify cell neighbourhoods. Building on this foundation, we apply this approach here to identify the cells themselves and explore the performance of the proposed method. Specifically, we initially apply this approach to the simplest case of identifying clique triangles, and then extend it to detect higher-order cliques. Throughout this process, we analyze the bottlenecks and key features of each classic algorithm, ultimately developing alternative methods that combine elements of traditional approaches \cite{tomita2003efficient,tomita2006worst} with the node neighbourhood strategy. More specifically, among the various algorithms for finding triangles \cite{al2018triangle, hu2014efficient, arifuzzaman2015space, arifuzzaman2015space, listing2016efficient, uddin2019finding}, we drew inspiration from those that iterate over edges. While these algorithms are effective for detecting low-dimensional cliques, they are limited in scope. For computing higher-order cliques in Vietoris-Rips complexes, the authors of \cite{boissonnat2014simplex} propose a simplified expansion that utilizes the simplicial tree method. However, this approach faces significant memory challenges, especially in dense networks, as the simplicial mapping must be stored. In response, we propose algorithms that integrate both methods, alongside the features of the aforementioned algorithms, and introduce an alternative set-theoretical approach inspired by our recent developments\cite{de2023efficient}. Finally, we benchmarked our proposed algorithm against the classic methods available in the literature. Remarkably, we not only developed new algorithms for computing cliques but also superseded the performance of existing methods found in the literature. Our findings represent a significant advancement in the computation of higher-order structures within complex networks.

The remainder of this article is organised as follows.
In Section~\ref{sec:network_backgroup}, we briefly recall the foundational concepts of network theory and the Vietoris-Rips complex. Then, in Section~\ref{sec:proposed_algorithms}, we develop novel algorithms, using a set-theoretical approach, for identifying cliques. Subsequently, in Section~\ref{sec:benchmark}, we benchmark our algorithms and compare their performance with classic implementations available in the literature. Finally, in Section~\ref{sec:conclusion}, we present our conclusions and propose future perspectives on this work.

\section{Basics of network theory}
\label{sec:network_backgroup}
For the sake of clarity, we first outline the essentials of network theory and the Vietoris-Rips complex. For further details, however, we refer the reader to the standard literature on network theory~\cite{west2001introduction,bondy2008graph}  and the Vietoris-Rips complex~\cite{zomorodian2005topology,zomorodian2010fast,edelsbrunner2022computational}.

\subsection{Undirected graph}
An \textit{undirected simple graph} $G=(V,E)$ is defined by a finite \textit{set of nodes (vertices)} $V$ and a \textit{set of edges} $E=\{\{x,y\}\,|\, x,y\in V, x\neq y\}$.
Let $x,y \in V$. We say that $y$ is a \textit{neighbour} of $x$ if the edge $(x,y) \in E$. The \emph{neighbourhood} of node $x$ is defined by a set of nodes that are connected to $x$ via an edge in $E$, and we denote \md{it} by $\pi(x)$. Formally,
\begin{equation}\label{eq:node_neighborhood}
    \pi(x)=\{y\in V\,|\, \{x,y\} \in E \}.
\end{equation}
A \textit{subgraph} of $G$ is a pair $G'=(V',E')$ such that $V'\subseteq V$ and $E'\subseteq E$, and that for all $\{x, y\} \in E'$  we have $x, y \in V'$. A \textit{complete subgraph} in $G$ is a subgraph such that all nodes are connected by an edge. The complete subgraphs with $|V'|=d$ nodes are called $d$-cliques (or cliques of dimension $d+1$.) We say that $G$ is \textit{weighted} when we associate a (positive) function $f:E\rightarrow \mathbb{R}$ to the set of edges, and $w=f(e)$ is the \textit{weight} of the edge $e$. 

\subsection{Vietoris-Rips Complex and Cliques}
Given $d+1$ nodes in a graph $G$, the \textit{$d$-simplex} associated with these nodes is the $d$-dimensional polytope that their convex hull forms. Hence, a $0$-simplex is a given node, a $1$-simplex is an edge connecting two nodes, a $2$-simplex is a filed triangle, a $3$-simplex is a filed tetrahedron, and so on. 
In this context, a \textit{face} of the $d$-simplex formed by $(d+1)$ nodes of $G$ is any simplex obtained by considering a subset of these nodes. Hence, any $d$-simplex possesses $\binom{d+1}{i+1}$ faces of dimension $i$. When $\sigma'$ is a face of $\sigma$ such that $|\sigma'|=|\sigma|-1$, $\sigma'$ is said to belong to the \textit{boundary} of $\sigma$, denoted by $\partial(\sigma)$.

A \textit{simplicial complex} $C$ is a set of simplices endowed with the hereditary property: given an element $\sigma \in C$ if $\sigma'$ is a face of $\sigma$, then $\sigma' \in C$. If $d$ is the maximum dimension of the simplices in $C$, then $C$ is said to be $d$-dimensional and it is called a simplicial $d$-complex. The underlying graph of $C$ is the graph induced by the $1$-simplices (edges) in $C$, i.e., the graph whose node-set is the set of $0$-simplices and whose edge set is the set of $1$-simplices in $C$. If the maximum dimension of $d$ of $C$ is $0$, then $C$ has only $0$-simplices ($1$-cliques or nodes), the case $d=1$ implies that $C$ has $0$ and $1$-simplices ($2$-cliques or edges), when $d=2$ it has up to $2$-simplices ($3$-cliques or triangles), and so on. Consequently, a $d$-simplex $\sigma$ can be induced by a $(d+1)$-clique made of $\binom{d+1}{2}$ edges.
In other words, one can associate a network with a simplicial complex by considering a \emph{clique complex}: each clique (complete subgraph) of $(d+1)$ vertices in the network is seen as a $d$-simplex of a simplicial complex.

Given a set of vertices $V \subseteq \mathbb{R}^m$ then let's consider, for any $x$ in $V$, the open ball centered in $x$ and of radius $\varepsilon >0$ (cutoff parameter): $B_{\varepsilon}(x) = \{ y\in \mathbb{R}^m \mid \text{dist}(x,y) \leq \varepsilon \}$, where $\text{dist}(x,y)$ is the Euclidean distance between two points $x,y\in \mathbb{R}^m$. Then, a \textit{Vietoris-Rips (VR) complex} $C$ is defined as a set of simplices whose balls only require pairwise intersection to be non-empty, that is,
\begin{equation}
\label{eq:vr}
\sigma = \{ x_i\}_i \in C \iff B_{\varepsilon}(x_i) \cap B_{\varepsilon}(x_j) \ne \emptyset, \forall x_i,x_j \in \sigma .
\end{equation}
This is equivalent to a clique complex of a graph with $V$ as vertex set in which every edge $(x_i, x_j)$ satisfies $d(x_i, x_j) \leq 2 \varepsilon$.

\begin{figure}[!htb]
    \centering
    \includegraphics[width=0.9\linewidth]{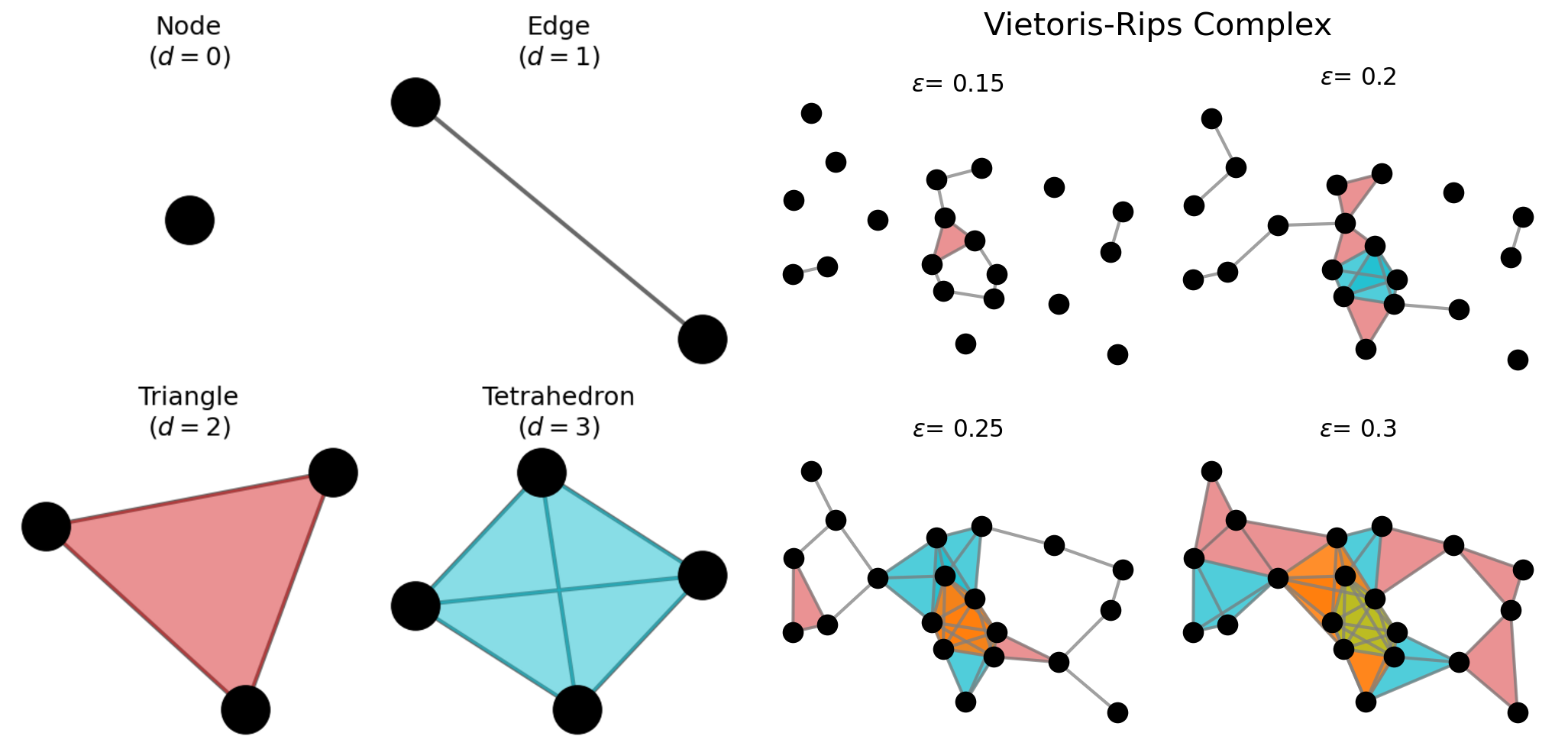}
    \caption{Example of cliques of dimensions $d \in \{0,1,2,3\}$, and the construction of a Vietories-Rips complex for different distance parameters \md{$\varepsilon$}. Higher dimensional cliques appear as the cutoff distance increases. For instance, tetrahedron (blue faces), pentahedron (orange faces) and hexahedron (green faces).}
    \label{fig:Cliques_VR}
\end{figure}

A \textit{finite filtration} that constructs a Vietoris-Rips complex generates a nested sequence of subcomplexes
\begin{equation} \label{eq.subcomplexes}  
\emptyset \subseteq C_{\varepsilon_0} \subseteq C_{\varepsilon_1} \subseteq \ldots  \subseteq C_{\varepsilon_n} \subseteq C,
\end{equation}
where $\varepsilon_0\leq \varepsilon_1\leq \hdots \leq \varepsilon_n$ and $C_\varepsilon=\{\sigma=\{x_i\}_i \,|\,\text{dist}(x_i,x_j)\leq 2\varepsilon, \forall x_i,x_j \in \sigma\}$. This leads to a weighted simplicial complex which is parameterised by the cutoff parameter $\varepsilon\geq 0$. Figure~\ref{fig:Cliques_VR} shows examples of cliques \md{of} different dimensions. This concept is the key point to define a nested chain of simplicial complexes that evolve according to the edges' weights.

\section{Proposed Set-theoretical based algorithms}\label{sec:proposed_algorithms}
Despite the efficient and natural construction of a Vietoris-Rips complex \cite{zomorodian2010fast}, the process of identifying cliques remains a significant challenge \cite{dieu1986average, fellows2009clique}. To address this issue, a wide range of algorithms have been proposed, spanning approaches from computing triangles \cite{chiba1985arboricity,schank2005finding,hu2014efficient,arifuzzaman2015space,listing2016efficient,uddin2019finding} to the computation of more complex structures \cite{bron1973algorithm,tomita2003efficient,zhang2005genome,tomita2006worst,boissonnat2014simplex,baudin2023faster}. Yet, each algorithm encounters its own set of challenges, such as memory leaks, computation time \cite{bron1973algorithm,zhang2005genome} or difficulties in implementation leading to solutions that are not scalable \cite{lagraa2016efficient}. To address these issues, we employ set-theoretical approaches~\cite{enderton1977elements,jech2003set,zimmermann1985applications,pawlak2002rough}, which offer alternative optimization strategies for tackling computational problems in computer science. Indeed, it facilitates efficient combinatorial computations and allows us to leverage well-optimized standard libraries for set-theory operations, enabling high-performance implementations. As an example of a recent successful application, we cite our efficient computation of Forman-Ricci curvature of simplicial complexes~\cite{Barros_de_Souza_2021}.

\begin{figure}
    \centering
    \includegraphics[width=0.9\linewidth]{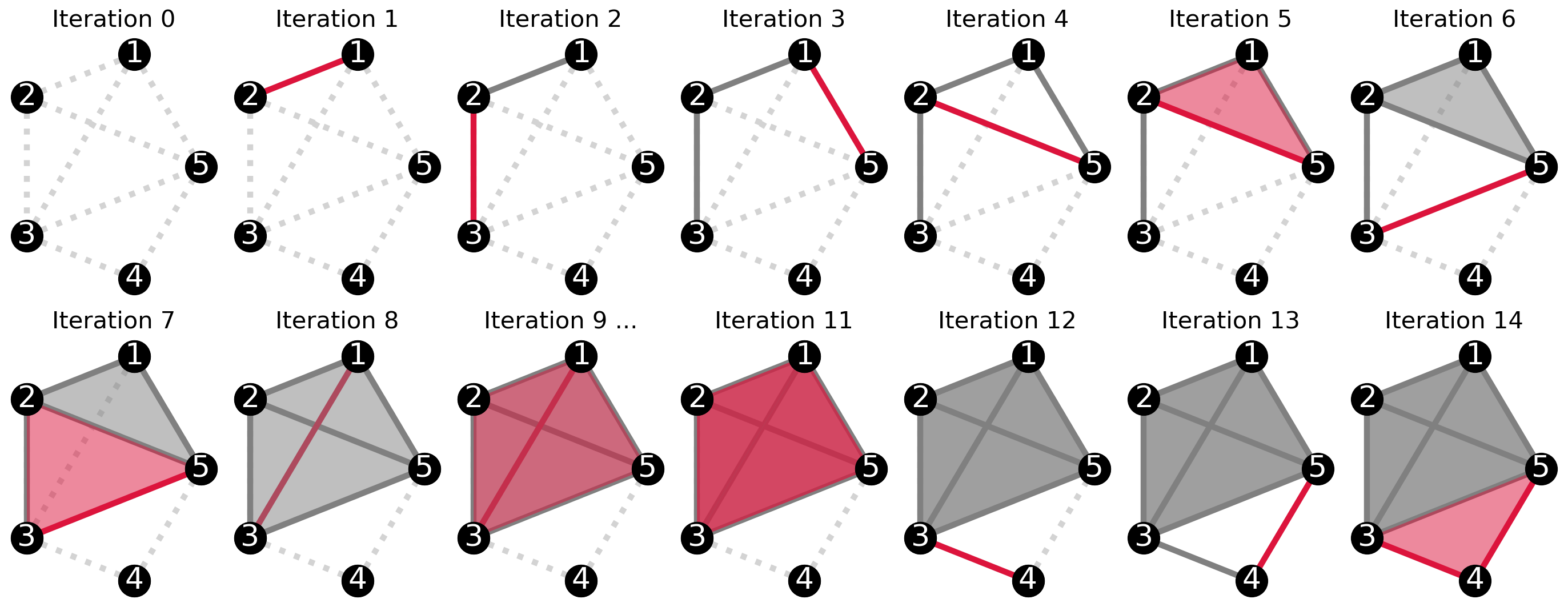}
    \caption{Finding cliques in a simple network as the neighbourhood is updated. The idea of the algorithms relies on the extension of the iterated edges to higher-order simplices. In this case, the ordering of edges is $(1,2),(2,3),(1,5),(3,5),(1,3),(3,4),(4,5)$. One new edge can give rise to several cliques, for instance, edge $(1,3)$ at iteration 8, which provides clique expansions at iterations $9$, $10$ (two triangles - not visible) and $11$ (one tetrahedron).}
    \label{fig:iteration_process}
\end{figure}

To advance, we systematically analyzed state-of-the-art algorithms, identifying both their advantages and critical bottlenecks. This allowed us to pinpoint key features that can be effectively scaled for higher-order cliques. Despite their varying methodologies, most existing algorithms share a fundamental limitation: a static view of node neighbourhoods. This leads to excessive clique searches, redundant data removal, and memory inefficiencies. To address this, we propose novel algorithms that dynamically update node adjacency based on filtration values, enabling efficient processing of evolving networks. Drawing inspiration from the simplicial tree approach, where simplicial paths are stored, we instead maintain local node neighbourhood information as cliques are updated. Our approach, rooted in set-theoretic discretization, builds upon our previous work ~\cite{Barros_de_Souza_2021}. By recursively extending node neighbourhoods to form complete subgraphs, we ensure real-time neighbourhood updates, significantly improving computational efficiency. This dynamic approach not only mitigates memory leaks but also optimizes clique identification, making it highly scalable for complex network structures. The visualization of this extension and its impact on clique detection is explored in detail in Figure~\ref{fig:iteration_process}. 

Moreover, a crucial factor for efficiency is identifying cliques at the instant when the neighbourhood is accessed. Building on this principle, we introduce two distinct methods to describe neighbourhood structures:

\begin{enumerate}
    \item[]\md{Method 1:} \md{Directly} from the neighboring nodes of the clique's boundary;
    \item[]\md{Method 2:} \md{Incrementing the node neighbourhood according to the clique's dimension within the classic node neighbourhood approach}.
\end{enumerate}
To \md{implement} these \md{methods}, we need to define separately the extensions of neighborhood definition for each approach. \md{Method $1$} will be called \textit{Dynamic neighbourhood boundary approach} \md{and method $2$} will be referred to as \textit{Dynamic Multi-layer Node-Neighbourhood approach}.
\md{In both cases}, the starting point for implementation for triangles \md{is described in} \md{A}lgorithm~\ref{alg:triangles_new_alg}\md{; see details in Appendix~\ref{sec:proprosed_algorithms_codes}}. For higher-order clique\md{s (}dimensions $d_{\max}> 2$), the approach \md{will be} split between the methods $1$ and $2$, which we \md{explain next}.
\subsection{Evolving Neighbourhood Boundary Approach}\label{sec:DBN_method}

In our previous work \cite{Barros_de_Souza_2021}, we demonstrated that for a $d$-simplex $\sigma$ in a simplicial complex $C$ (with dimension greater than $d$), a $(d+1)$-simplex containing $\sigma$ as a face can be identified using the intersection of the boundary's neighbourhood in the underlying graph. For instance, if $\sigma$ is a $2$-simplex with vertices $\{x_1, x_2, x_3\}$, and all these nodes share a common neighbour $x_4$ in the graph, then $\{x_1, x_2, x_3, x_4\}$ naturally forms a $3$-simplex in $C$, with $\sigma$ as one of its four faces.

Following this idea, we expand the neighbourhood of cliques based on how \md{their} boundaries behave, that is, by studying the neighbourhood of each boundary. This information is stored dynamically as the network evolves \md{at} each step of edge iteration. As an example, for {identifying} tetrahedr\md{a}, we need the information of the previous sequence that generated the triangles, that is, we need the information of the previous triangles as the starting point so that the subsequent list of tetrahedr\md{a} can be found from the neighbourhood intersection. \md{The same idea underpins} the forward method~\cite{schank2005finding}, where the triangles are found as an extension of edge iterations. 

The optimization challenge lies in providing extensions while avoiding repetitions and storing redundant information. To this end, we update the neighbourhood of the \md{clique} boundaries \md{rather} than dealing with the neighbourhood of nodes as in the initial approach. More precisely, let $\triangle=\{x,y,z\}$ be a triangle provided by a clique extension\md{, that is,} $\triangle=\{x,y\}\cup\{z\}$, where $\{x,y\}$ was the original edge iterated and $z$ is the node \md{that} extends the edge to the triangle. For each element $e \in \triangle$, we update the neighbourhood of tetrahedr\md{a} boundaries, $\pi_{\tetrahedron}(e):=\pi_{\tetrahedron}(e)\cup\{\triangle\setminus e\}$, and then we iterate the intersection of the triangles boundaries for finding the tetrahedrons that extend $\triangle$. This idea is \md{formalised} in \md{A}lgorithm~\ref{alg:tetra_dn}. Subsequently, for each tetrahedron $\tetrahedron$, the extension for {identifying} pentahedr\md{a follows} the same \md{strategy} for updating the neighbourhood and finding the set of nodes that \md{allows this extension from the initial} tetrahedron. Explicitly, for each $\triangle \in \,\tetrahedron$, we update $\pi_{4}(\triangle):=\pi_{4}(\triangle)\cup\{\tetrahedron\setminus \triangle\}$ and iterate each node of the intersection of neighbours per boundary in $\tetrahedron$\md{, which extends} the tetrahedron to pentahedr\md{a}. \md{A}lgorithm~\ref{alg:penta_dn} describes these dynamics. \md{This} idea can \md{naturally} be \md{extended to} hexahedr\md{a} (\md{see A}lgorithm~\ref{alg:hexa_dn}) and so forth. \md{Then, a systematic extension to} cliques up to dimension $d_{\max}$ \md{is achieved via A}lgorithm~\ref{alg:RDBNA}.
{To elucidate the core of our algorithmic insight}, Table \ref{tab:cliques_boundaries_example} provides the information of cliques found, as well as, the boundary of each clique per iteration step, for the example {depicted} in Figure \ref{fig:iteration_process}. {Moreover,} Table \ref{tab:boundary-set} provides the evolution of the set of boundary neighbours per iteration.
\subsection{Evolving Multi-Layer Node-Neighbourhood Approach}\label{sec:DMNN_method}
\md{Here, we adopt a slightly different approach.} {Specifically,} instead of storing information about clique neighbourhoods from their boundary, we map the higher-order cliques' information based on the neighbouring nodes whose intersection provides the clique extensions to a higher-order clique (as performed \md{in A}lgorithm~\ref{alg:triangles_new_alg}). However, we split the neighbourhood information according to their dimension, which is provided as layers per dimension of neighbourhood information\md{. C}onsequently, the neighbourhood is updated separately, \md{hence} avoiding clique repetition. 

To exemplify, consider the case of identifying tetrahedra $\tetrahedron$. First assume we are given the algorithmic steps to finding triangles $\triangle=(x,y,z)$ (originally extended from an iterated edge, $e=(x,y)$). Moreover, assume the neighbourhood for triangles is updated according to Algorithm~\ref{alg:triangles_new_alg}. To update the neighbourhood layer for tetrahedr\md{a} ($\pi_{\tetrahedron}$) we iterate all the nodes $c \in e$ and update by adding the last element in for generating the triangle, i.e., $\pi_{\tetrahedron}:=\pi_{\tetrahedron}\cup \{z\}$. The intersection of all neighbourhood nodes in $\triangle$ \md{with}in the $\pi_{\tetrahedron}$ will provide the set of nodes that extends $\triangle$ to a tetrahedron\md{; A}lgorithm \ref{alg:tetra_nb} details this process.

\md{T}he extensions for identifying pentahedr\md{on} (cliques of dimension $4$) from a given tetrahedron $\tetrahedron=(x,y,z,w)$, can naturally be achieved by updating the neighbourhood layer for the pentahedron. Specifically, for each node $n \in \triangle$, we refine its neighbourhood by update $\pi_4(n):=\pi_4(n)\cup \{w\}$. Further, we compute the intersection of $\pi_4(n)$ \md{for all $n$ in $\tetrahedron$} to find all nodes that extend to pentahedr\md{a; this corresponds to A}lgorithm~\ref{alg:penta_nb}. Analogously, we can find extensions for hexahedr\md{on, as described in A}lgorithm~\ref{alg:hexa_nb}. In a recurrent way, we can extend \md{this approach to} cliques up to dimension $d_{\max}$ \md{as described in A}lgorithm~\ref{alg:RDNNA}. To further elucidate our algorithmic strategy, we provide in Table \ref{tab:node-neigh} the identified cliques in each iteration for the example depicted in Figure \ref{fig:iteration_process}.

%
\section{Benchmark Report}
\label{sec:benchmark}
To evaluate the efficiency of the proposed algorithms, whose pseudocodes are provided in Appendix~\ref{sec:proprosed_algorithms_codes}, we implemented them in Python~\cite{python} and compared their performance with two widely used algorithms for finding cliques, namely, \texttt{NetworkX} \cite{hagberg2020networkx}, which uses Zhang's algorithm~\cite{zhang2005genome} and \texttt{Gudhi}~\cite{maria2014gudhi}, which employs the simplicial tree method~\cite{boissonnat2014simplex}. Noteworthy, for our test, we adapted the \texttt{NetworkX} code to allow the search for cliques to be restricted by a specified maximum dimension detailed code available at \cite{souza2025benchmark}). The benchmarks were performed on a Dell XPS 15-7590 laptop equipped with an Intel Core i7-9750H CPU and 32GB of RAM, running Ubuntu 20.04.4 LTS. The analysis was divided into two categories: applications using synthetic data and those involving point-cloud data.

\subsection{Performance in  VR-complex filtrations generated from complete graphs}
\begin{figure}[!htb]
    \centering
    \includegraphics[width=0.9\linewidth]{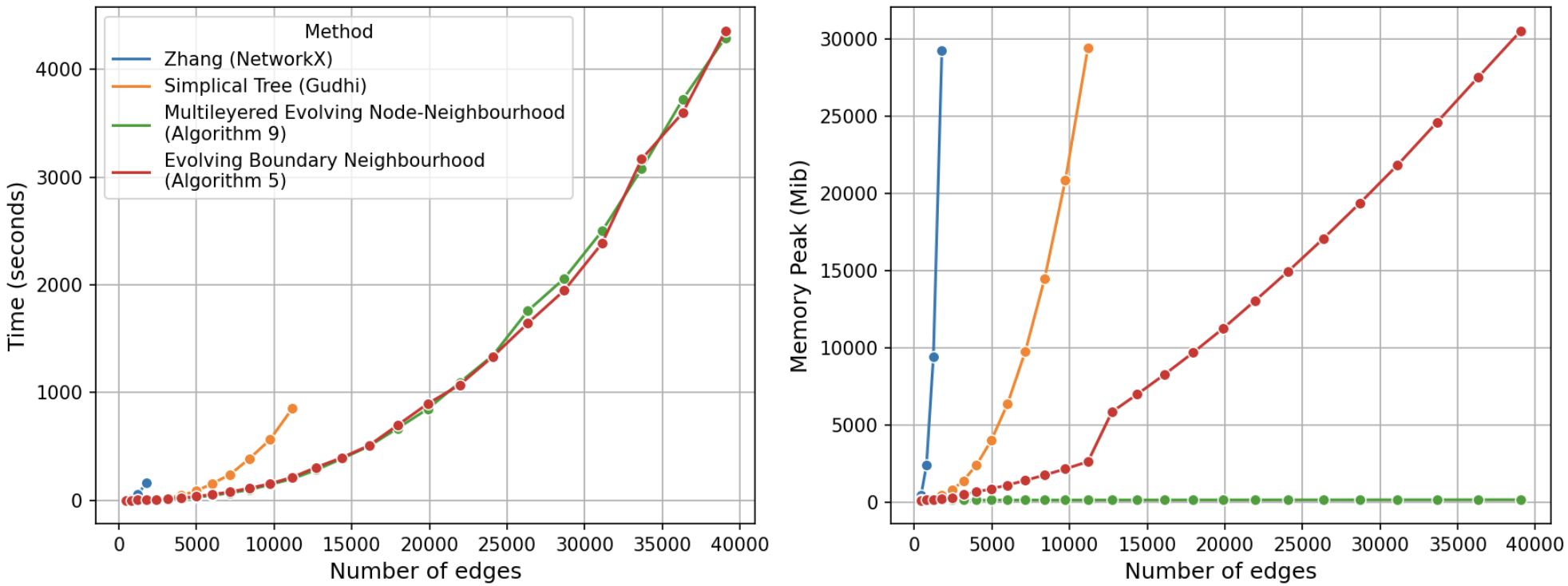}
    \caption{Benchmark comparison between algorithms for finding cliques in VR-complexes from complete networks ($d_{\max}$=4).}
    \label{fig:Benchmark_complete_graphs}
\end{figure}

As an initial step in our testing process, we generated complete graphs with node counts ranging from $30$ to $150$ with increments of $10$ between the two bounds. We ran the algorithms under stress conditions with complete (weighted) graphs, where the filtration process runs until the complete graph is obtained. We limited the search for finding cliques to the maximum dimension $d_{\max}=4$ (up to pentahedron). As depicted in Figure~\ref{fig:Benchmark_complete_graphs}, the \texttt{NetworkX} algorithm displays the poorest performance in clique identification, primarily due to memory leaks, which is a significant bottleneck. In contrast, the \texttt{Gudhi}'s clique algorithm exhibits a notable improvement for complete weighted graphs with up to $n=10,000$ nodes (\textit{i.e.,  $\binom{n}{2}$} edges). However, beyond this threshold, memory consumption increases sharply, also negatively impacting processing time. Remarkably, our algorithms outperform both in terms of time and memory efficiency, enabling significantly faster computations and broadening the potential for computing cliques in denser networks. This fact demonstrates the efficiency of our proposed algorithm on computing cliques along large cutoff distances in small dense graphs.

\subsection{Performance in VR-complexes generated from  point cloud data}
\begin{figure}
    \centering
    \includegraphics[width=0.9\linewidth]{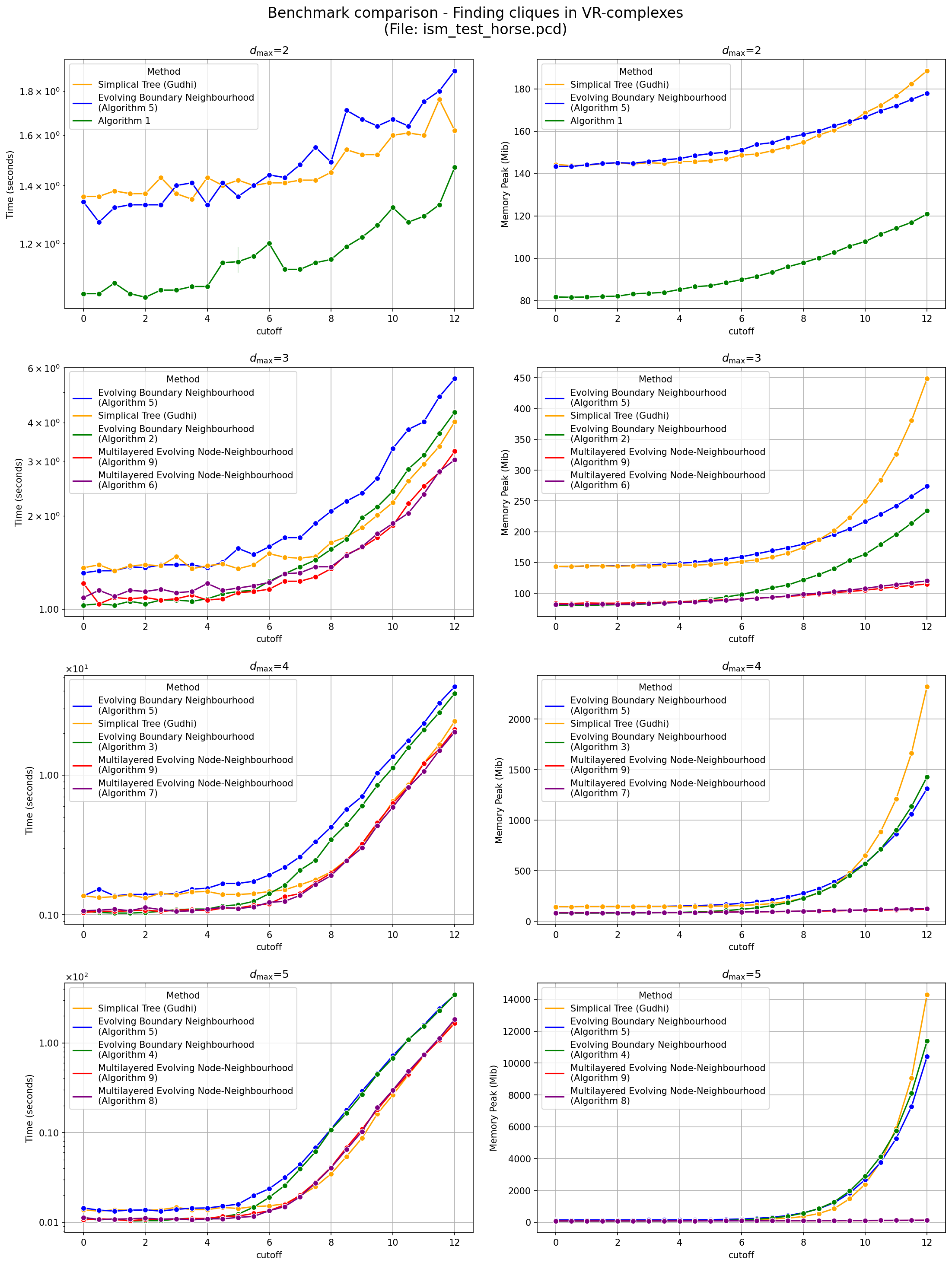}
    \caption{Benchmark comparison between algorithms for finding cliques in VR-complexes from point cloud data ((\textbf{ism\_test\_horse.pcd})).}
    \label{fig:Benchmark_horse}
\end{figure}

\begin{figure}
    \centering
    \includegraphics[width=0.9\linewidth]{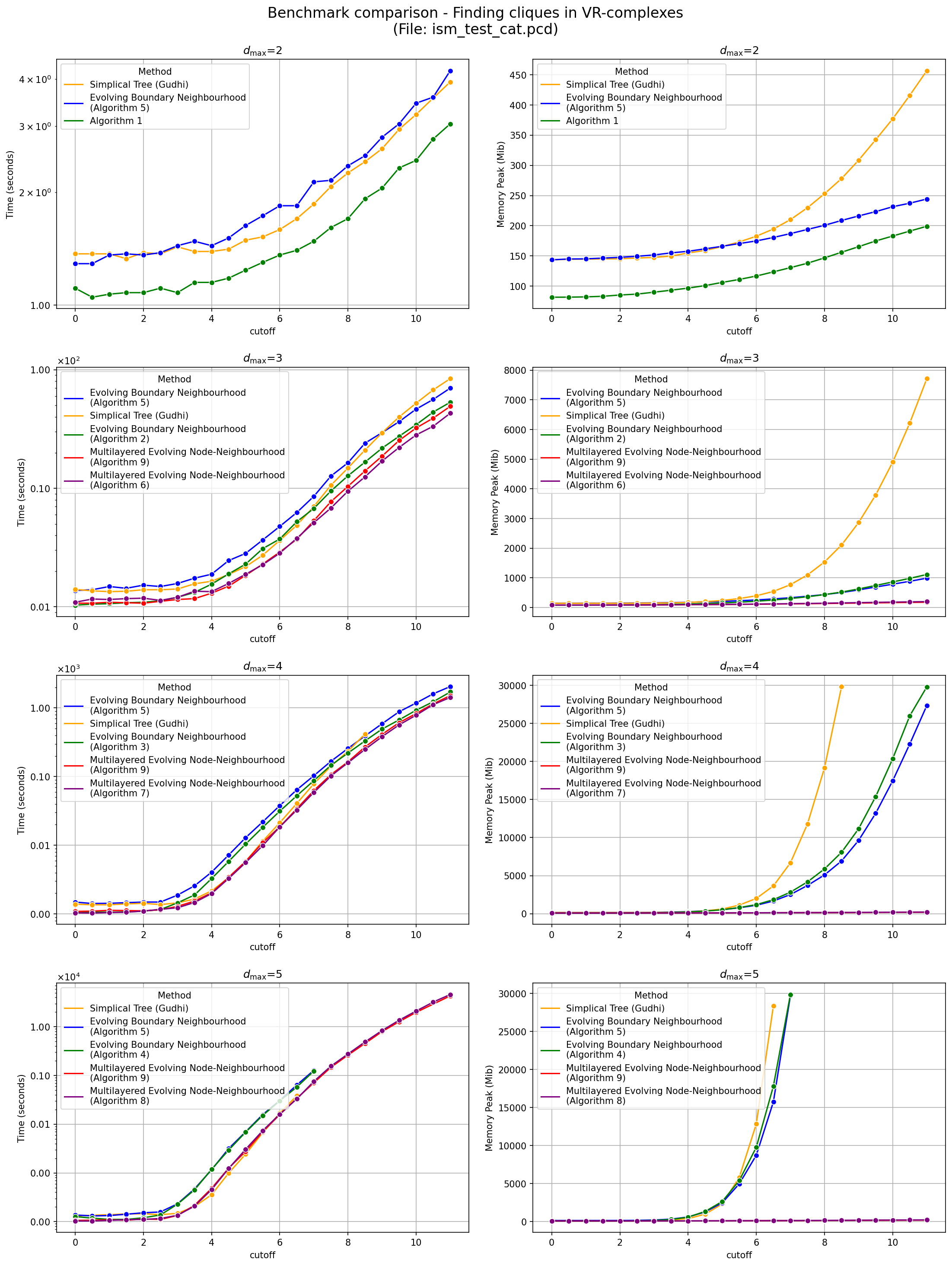}
    \caption{Benchmark comparison between algorithms for finding cliques in VR-complexes from point cloud data ( \textbf{ism\_test\_cat.pcd}).}
    \label{fig:Benchmark_cat}
\end{figure}

To evaluate our algorithms using point cloud data, we utilized the files \textbf{ism\_test\_horse.pcd} and \textbf{ism\_test\_cat.pcd}, each containing 3,400 nodes (see the Point Cloud Libary repository \cite{PCLData}). For this analysis, we excluded the \texttt{NetworkX} algorithm due to performance issues. We repeated the tests as a function of the cutoff distance between nodes, for the two sets of data and for the clique dimensions $2$ (triangles), $3$ (tetrahedrons), $4$ (pentahedrons), and $5$ (hexahedrons). In all cases, our proposed algorithms outperformed the memory efficiency.  The results are shown in Fig. \ref{fig:Benchmark_horse} and Fig. \ref{fig:Benchmark_cat}. Strikingly, the Multilayered Node-Neighbourhood algorithms, particularly Algorithm $9$, demonstrated a significant improvement in preserving neighbourhood information, resulting in the lowest memory leak when compared to the other algorithms analyzed. The drawback of this approach is the inflexibility in choosing the edge ordering, which is required, for instance, in persistence diagram approaches. In contrast, the evolving boundary neighbourhood algorithms (such as Algorithm $5$) allow for the selection of an edge ordering; however, time and memory efficiency are compromised for high-dimensional cliques and higher cutoff values. This is due to the large amount of boundary information that must be stored as the number of cliques increases under these conditions. However, the memory increase is not as significant as with the \texttt{Gudhi} algorithm. Overall, our proposed algorithms can be chosen based on the features required for constructing VR complexes and hardware limitations.
\subsection{Performance in unweighted Real-World Networks} 
\begin{figure}
    \centering
    \includegraphics[width=0.9\linewidth]{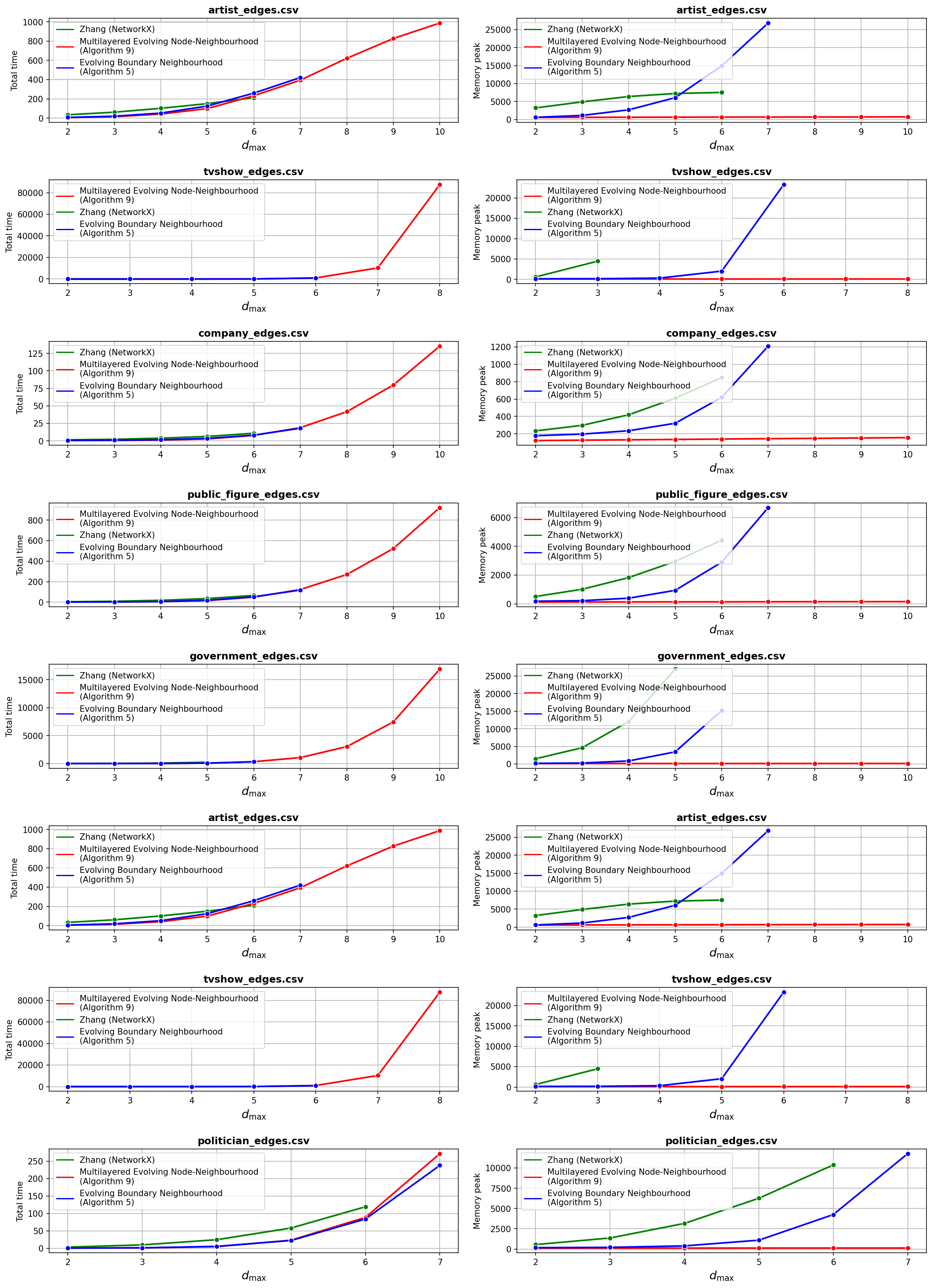}
    \caption{Benchmark comparison between algorithms for finding cliques in VR-complexes from real data unweighted networks (Facebook data from SNAP repository \cite{rozemberczki2019gemsec}).}
    \label{fig:Benchmark_unweighted}
\end{figure}

To strengthen the performance analysis of our proposed algorithms, we extended the benchmark to unweighted real-world networks. Given Gudhi’s limitation to Vietoris-Rips complexes, we compared our algorithms with a modified version of Zhang’s algorithm (NetworkX). For this, we utilized Facebook data from the SNAP repository \cite{rozemberczki2019gemsec}, with results shown in Figure \ref{fig:Benchmark_unweighted}. We evaluated algorithm efficiency based on the maximum clique dimension. Consistent with prior results, our algorithms outperformed others in memory management. Notably, Algorithm \ref{alg:RDNNA} achieved the best memory efficiency, ensuring scalability to higher clique dimensions without compromising performance. This advantage arises from storing neighbourhood information per node rather than using Zhang’s method or the boundary storage approach in Algorithm \ref{alg:RDBNA}. Table \ref{tab:methods_comparison} summarizes the benchmark, highlighting the strengths and limitations of each method.

\begin{table}[h!]
\centering
\begin{tabular}{|l|l|l|l|l|}
\hline
Method & Zhang & Simplicial Tree & Neigh. Boundary & Multi-Layered Node-Neighbourhood \\
\hline
Algorithm & NetworkX & Gudhi & Algorithms 2,3,4 and 5 & Algorithms 6,7,8 and 9 \\
\hline
Memory Efficient & No & No & No & Yes \\
\hline
Sorts Cliques by Weight & No & Yes & No & Yes \\
\hline
\end{tabular}
\caption{Summary of the methods bench-marked in our analysis.}
\label{tab:methods_comparison}
\end{table}

\newpage
\section{Conclusion}\label{sec:conclusion}
Network science which includes higher-order networks, is in pursuit of highly efficient algorithms to solve complex challenges across science, engineering, healthcare, and economics. At the core of these networks lies the concept of cliques, a fundamental structure that helps organize the inherent complexity. Yet, the identification of cliques remains an NP-complete problem, presenting a significant computational challenge. To advance, we begin by drawing inspiration from state-of-the-art algorithms for clique identification and analyzing their strengths and limitations. A key drawback of these algorithms is their static view of the node neighbourhood, which induces excessive searching for cliques and redundant information, resulting in inefficient runtime and memory usage. To overcome these challenges, we build on the strengths of our previous work~\cite{Barros_de_Souza_2021}, casting these algorithms in the context of set theory. Subsequently, we generalize the identification process to higher-order cliques (specifically cliques within Vietoris-Rips complexes) and improve it by discovering novel ways to explore the neighbourhoods of cliques. In particular, we present two key paradigms, which we denote as \textit{evolving neighbourhood boundary} and \textit{evolving multi-layer node neighbourhood}. Further, we leverage well-optimized standard libraries for set-theory operations, enabling Python-based high-performance implementations. Finally, we benchmark and demonstrate that our proposed algorithms provide significant memory efficiency in the context of cliques in dense networks, thus superseding state-of-the-art algorithms. Consequently, we believe that the outcomes of our work will impact big data and complex network science, in particular in the context of topological and geometrical data analysis. For instance, our proposed algorithms can be employed as alternatives for computing the \textit{Euler Characteristics} and their associated transitions (e.g. via topological phase transitions) in complex systems (e.g. Brain data). Furthermore, the correct choice in the clique's method computation has the potential use of identifying biological markers (e.g. in personalized medicine) and clustering information from a geometric viewpoint (e.g., by using \textit{Forman-Ricci curvature} computations.)
\section*{Acknowledgments}
\textcolor{black}{This research is supported by the grant PID2023-146683OB-100 funded by MICIU/AEI /10.13039/501100011033 and by ERDF, EU. Additionally, it is supported by Ikerbasque Foundation and the Basque Government through the BERC 2022-2025 program and by the Ministry of Science and Innovation: BCAM Severo Ochoa accreditation \newline CEX2021-001142-S / MICIU / AEI / 10.13039/501100011033. Moreover, the authors acknowledge the financial support received from BCAM-IKUR, funded by the Basque Government by the IKUR Strategy and by the European Union NextGenerationEU/PRTR. We also acknowledge the support of ONBODY no. KK-2023/00070 funded by the Basque Government through ELKARTEK Programme.}
%
%
%
\appendix
\section{\md{Pseudo-codes}}
\label{sec:appendix}
In this section, we provide the pseudocodes derived from our theoretical and computational analysis and for completeness the classic methods that inspired our work.
\subsection{\md{Proposed Algorithms}}
\label{sec:proprosed_algorithms_codes}
\begin{algorithm}[H]
\caption{Find Triangles (Evolving Neighborhood Approach)}\label{alg:triangles_new_alg}
\begin{algorithmic}[1]
\Require List of edges, $E$
\State $\pi=\emptyset$
\For{$e=(x,y) \in E$ }
    
    \State{$\pi(x):=\pi(x)\cup\{y\}$}
    \State{$\pi(y):=\pi(y)\cup\{x\}$}
   \State{$Z=\pi(x)\cap\pi(y)$}
   \For{$z \in Z$}
   \State{\textbf{print} $\{x,y,z\}$}
   \EndFor
\EndFor
\end{algorithmic}
\end{algorithm}

\begin{algorithm}[H]
\caption{Find Triangles and Tetrahedrons (Evolving Boundary Neighborhood Approach)}\label{alg:tetra_dn}
\begin{algorithmic}[1]
\Require List of (sorted) edges, $E$
\State {$\pi_{\triangle}=\emptyset,\pi_{\tetrahedron}=\emptyset$}
\For{$e=(x,y) \in E$ }
    
    \State{$\pi_{\triangle}(x):=\pi_{\triangle}(x)\cup\{y\}$}
    \State{$\pi_{\triangle}(y):=\pi_{\triangle}(y)\cup\{x\}$}
   \State{$Z=\pi_{\triangle}(x)\cap\pi_{\triangle}(y)$}
   \For{$z \in Z$}
    \State{$\triangle=\{x,y,z\}$}
   \State{\textbf{print} $\triangle$}
   \State{$\partial(\triangle)=\{\{x,y\},\{x,z\},\{y,z\}\}$}
   \State{$\pi_{\tetrahedron}(\{x,y\}):=\pi_{\tetrahedron}(\{x,y\})\cup \{z\}$}
   \State{$\pi_{\tetrahedron}(\{x,z\}):=\pi_{\tetrahedron}(\{x,z\})\cup \{y\}$}
   \State{$\pi_{\tetrahedron}(\{y,z\}):=\pi_{\tetrahedron}(\{y,z\})\cup \{x\}$}
   \State{$W=\pi_{\tetrahedron}(\{x,y\})\cap\pi_{\tetrahedron}(\{x,z\})\cap \pi_{\tetrahedron}(\{y,z\})$}
   \For{$w \in W$}
   \State{$\tetrahedron=\{x,y,z,w\}$}
   \State{\textbf{print} $\tetrahedron$}
   \EndFor
   \EndFor
\EndFor
\State 
\end{algorithmic}
\end{algorithm}

\begin{algorithm}[H]
\caption{Find up to pentahedrons (Evolving Boundary Neighborhood Approach)}\label{alg:penta_dn}
\begin{algorithmic}[1]
\Require List of (sorted) edges, $E$
\State {$\pi_2=\emptyset,\pi_3=\emptyset,\pi_4=\emptyset$}
\For{$e=(x,y) \in E$ }
    
    \State{$\pi_{2}(x):=\pi_{2}(x)\cup\{y\}$}
    \State{$\pi_{2}(y):=\pi_{2}(y)\cup\{x\}$}
   \State{$Z=\pi_{2}(x)\cap\pi_{2}(y)$}
   \For{$z \in Z$}
    \State{$\tau_2=\{x,y,z\}$}
   \State{\textbf{print} $\tau_2$}
   \State{$\partial(\tau_2)=\{\{x,y\},\{x,z\},\{y,z\}\}$}
   \State{$\pi_{3}(\{x,y\}):=\pi_{3}(\{x,y\})\cup \{z\}$}
   \State{$\pi_{3}(\{x,z\}):=\pi_{3}(\{x,z\})\cup \{y\}$}
   \State{$\pi_{3}(\{y,z\}):=\pi_{3}(\{y,z\})\cup \{x\}$}
   \State{$W=\pi_{3}(\{x,y\})\cap\pi_{3}(\{x,z\})\cap \pi_{3}(\{y,z\})$}
   \For{$w \in W$}
   \State{$\tau_3=\{x,y,z,w\}$}
   \State{\textbf{print} $\tau_3$}
    \State{$\partial(\tau_3)=\{\{x,y,z\},\{x,y,w\},\{x,z,w\},\{y,z,w\}\}$}
    \State{$\pi_4(\{x,y,z\}):=\pi_4(\{x,y,z\})\cup \{w\}$}
    \State{$\pi_4(\{x,y,w\}):=\pi_4(\{x,y,w\})\cup \{z\}$}
    \State{$\pi_4(\{x,z,w\}):=\pi_4(\{x,z,w\})\cup \{y\}$}
    \State{$\pi_4(\{y,z,w\}):=\pi_4(\{y,z,w\})\cup \{x\}$}
    \State{$Q=\pi_4(\{x,y,z\})\cap \pi_4(\{x,y,w\})\cap \pi_4(\{x,z,w\})\cap \pi_4(\{y,z,w\})$}
    \For{$q \in Q$}
    \State{$\tau_4=\{x,y,z,w,q\}$}
    \State{\textbf{print} $\tau_4$}
    \EndFor
   \EndFor
   \EndFor
\EndFor
\end{algorithmic}
\end{algorithm}

\begin{algorithm}[H]
\caption{Find up to hexahedrons (Evolving Boundary Neighborhood Approach)}\label{alg:hexa_dn}
\begin{algorithmic}[1]
\Require List of (sorted) edges, $E$
\State {$\pi_2=\emptyset,\pi_3=\emptyset,\pi_4=\emptyset, \pi_5=\emptyset$}
\For{$e=(x,y) \in E$ }
    
    \State{$\pi_{2}(x):=\pi_{2}(x)\cup\{y\}$}
    \State{$\pi_{2}(y):=\pi_{2}(y)\cup\{x\}$}
   \State{$V=\pi_{2}(x)\cap\pi_{2}(y)$}
   \For{$z \in Z$}
    \State{$\tau_2=\{x,y,z\}$}
   \State{\textbf{print} $\tau_2$}
   \State{$\partial(\tau_2)=\{\{x,y\},\{x,z\},\{y,z\}\}$}
   \State{$\pi_{3}(\{x,y\}):=\pi_{3}(\{x,y\})\cup \{z\}$}
   \State{$\pi_{3}(\{x,z\}):=\pi_{3}(\{x,z\})\cup \{y\}$}
   \State{$\pi_{3}(\{y,z\}):=\pi_{3}(\{y,z\})\cup \{x\}$}
   \State{$W=\pi_{3}(\{x,y\})\cap\pi_{3}(\{x,z\})\cap \pi_{3}(\{y,z\})$}
   \For{$w \in W$}
   \State{$\tau_3=\{x,y,z,w\}$}
   \State{\textbf{print} $\tau_3$}
    \State{$\partial(\tau_3)=\{\{x,y,z\},\{x,y,w\},\{x,z,w\},\{y,z,w\}\}$}
    \State{$\pi_4(\{x,y,z\}):=\pi_4(\{x,y,z\})\cup \{w\}$}
    \State{$\pi_4(\{x,y,w\}):=\pi_4(\{x,y,w\})\cup \{z\}$}
    \State{$\pi_4(\{x,z,w\}):=\pi_4(\{x,z,w\})\cup \{y\}$}
    \State{$\pi_4(\{y,z,w\}):=\pi_4(\{y,z,w\})\cup \{x\}$}
    \State{$Q=\pi_4(\{x,y,z\})\cap \pi_4(\{x,y,w\})\cap \pi_4(\{x,z,w\})\cap \pi_4(\{y,z,w\})$}
    \For{$q \in Q$}
    \State{$\tau_4=\{x,y,z,w,q\}$}
    \State{\textbf{print} $\tau_4$}
    \State{$\partial(\tau_4)=\{\{x,y,z,w\},\{x,y,z,q\},\{x,y,w,q\},\{x,z,w,q\},\{y,z,w,q\}\}$}
    \State{$\pi_5(\{x,y,z,w\}):=\pi_5(\{x,y,z,w\})\cup\{q\}$}
    \State{$\pi_5(\{x,y,z,q\}):=\pi_5(\{x,y,z,q\})\cup\{w\}$}
    \State{$\pi_5(\{x,y,w,q\}):=\pi_5(\{x,y,w,q\})\cup\{z\}$}
    \State{$\pi_5(\{x,z,w,q\}):=\pi_5(\{x,z,w,q\})\cup\{y\}$}
    \State{$\pi_5(\{y,z,w,q\}):=\pi_5(\{y,z,w,q\})\cup\{x\}$}
    \State{$R=\pi_5(\{x,y,z,w\})\cap \pi_5(\{x,y,z,q\}) \cap \pi_5(\{x,y,w,q\})\cap \pi_5(\{x,z,w,q\})\cap \pi_5(\{y,z,w,q\})$}
    \For{$r \in R$}
    \State{$\tau_5=\{x,y,z,w,q,r\}$}
    \State{\textbf{print} $\tau_5$}
    \EndFor
    \EndFor
   \EndFor
   \EndFor
\EndFor
\end{algorithmic}
\end{algorithm}
\begin{algorithm}[H]
\caption{Compute Cliques (Recursive Evolving Boundary Neighbourhood Approach)}\label{alg:RDBNA}
\begin{algorithmic}[1]
    \Procedure{ComputeCliques}{$E$,$d_{\max}$}
        \State \textbf{Input:} $E$, the list of (sorted) edge, $d_{\max}$, the maximum clique dimension
        \State \textbf{Output:} All the cliques up to dimension $d_{\max}$
        \State \textbf{Initialization:}
        \State $\pi_d=\emptyset, \, \forall d \in \{2,\hdots,d_{\max}\}$
        \Statex

    \Procedure{AuxiliaryFunction}{$\sigma$}
        \State \textbf{Input:} $\sigma$, a clique
        \State \textbf{Output:} A sequence of cliques $\tau$ extended from $\sigma$
        \State $d=|\sigma|$
        \State \textbf{print} $\sigma$
        \State
        \If{$d<d_{\max}$}
        \For{$\gamma \in \partial(\sigma)$}
        \State $\pi_{d}(\gamma):=\pi_{d}(\gamma)\cup\{\sigma\setminus \gamma\}$
        \EndFor
        \State{$Z=\cap_{\gamma \in \partial(\sigma)} \pi_d(\gamma)$}
        \For{$z \in Z$}
        \State $\tau=\sigma \cup \{z\}$
        \State \textbf{compute} \textproc{AuxiliaryFunction}($\tau$)
        
        \EndFor
        \EndIf
    \EndProcedure

        \State
        \For{$e \in E$}
        \State \textbf{compute} \textproc{AuxiliaryFuncion}(e)
        \EndFor
        \State
    \EndProcedure
    
\end{algorithmic}
\end{algorithm}
\subsection{Evolving Multilayered Node-Neighbourhood Algorithms}
\begin{algorithm}[H]
\caption{Find Triangles and Tetrahedrons (Evolving Multilayered Node-Neighborhood Approach)}\label{alg:tetra_nb}
\begin{algorithmic}[1]
\Require List of edges (sorted in lexicographic order), $E$
\State {$\pi_{\triangle}=\emptyset,\pi_{\tetrahedron}=\emptyset$}
\For{$e=(x,y) \in E$ }
    
    \State{$\pi_{\triangle}(x):=\pi_{\triangle}(x)\cup\{y\}$}
    \State{$\pi_{\triangle}(y):=\pi_{\triangle}(y)\cup\{x\}$}
   \State{$Z=\pi_{\triangle}(x)\cap\pi_{\triangle}(y)$}
   \For{$z \in Z$}
    \State{$\triangle=\{x,y,z\}$}
   \State{\textbf{print} $\triangle$}
   \For{$c \in e$}
   \State{$\pi_{\tetrahedron}(c):=\pi_{\tetrahedron}(c)\cup \{z\}$}
   \EndFor
   \State{$W=\pi_{\tetrahedron}(x)\cap\pi_{\tetrahedron}(y)\cap \pi_{\tetrahedron}(z)$}
   \For{$w \in W$}
   \State{$\tetrahedron=\{x,y,z,w\}$}
   \State{\textbf{print} $\tetrahedron$}
   \EndFor
   \EndFor
\EndFor
\end{algorithmic}
\end{algorithm}

\begin{algorithm}[H]
\caption{Find up to Pentahedrons (Evolving Multilayered Node-Neighborhood Approach)}\label{alg:penta_nb}
\begin{algorithmic}[1]
\Require List of edges (sorted in lexicographic order), $E$
\State {$\pi_{2}=\emptyset,\pi_{3}=\emptyset,\pi_{4}=\emptyset$}
\For{$e=(x,y) \in E$ }
    
    \State{$\pi_{2}(x):=\pi_{2}(x)\cup\{y\}$}
    \State{$\pi_{2}(y):=\pi_{2}(y)\cup\{x\}$}
   \State{$Z=\pi_{2}(x)\cap\pi_{2}(y)$}
   \For{$z \in Z$}
    \State{$\tau_{2}=\{x,y,z\}$}
   \State{\textbf{print} $\tau_2$}
   \For{$n \in e$}
   \State{$\pi_{3}(n):=\pi_{3}(n)\cup \{z\}$}

   \EndFor
   \State{$W=\pi_{3}(x)\cap\pi_{3}(y)\cap \pi_{3}(z)$}
   \For{$w \in W$}
   \State{$\tau_4=\{x,y,z,w\}$}
   \State{\textbf{print} $\tau_3$}
   \For{$n \in \tau_2$}
   \State{$\pi_{4}(n):=\pi_{4}(n)\cup \{w\}$}
   \EndFor
   \State{$Q=\pi_{4}(x)\cap\pi_{4}(y)\cap \pi_{4}(z)\cap\pi_4(w)$}
 \For{$q \in Q$}
 \State{$\tau_4=\{x,y,z,w,q\}$}
 \State{\textbf{print} $\tau_4$}
 \EndFor
   \EndFor
   \EndFor
\EndFor
\end{algorithmic}
\end{algorithm}

\begin{algorithm}[H]
\caption{Find up to Hexahedrons (Evolving Multilayered Node-Neighborhood Approach)}\label{alg:hexa_nb}
\begin{algorithmic}[1]
\Require List of edges (sorted in lexicographic order), $E$
\State {$\pi_{2}=\emptyset,\pi_{3}=\emptyset,\pi_{4}=\emptyset,\pi_{5}=\emptyset$}
\For{$e=(x,y) \in E$ }
    
    \State{$\pi_{2}(x):=\pi_{2}(x)\cup\{y\}$}
    \State{$\pi_{2}(y):=\pi_{2}(y)\cup\{x\}$}
   \State{$Z=\pi_{2}(x)\cap\pi_{2}(y)$}
   \For{$z \in Z$}
    \State{$\tau_{2}=\{x,y,z\}$}
   \State{\textbf{print} $\tau_2$}
   \For{$n \in e$}
   \State{$\pi_{3}(n):=\pi_{3}(n)\cup (\{z\}$}
   \EndFor
   \State{$W=\pi_{3}(x)\cap\pi_{3}(y)\cap \pi_{3}(z)$}
   \For{$w \in W$}
   \State{$\tau_4=\{x,y,z,w\}$}
   \State{\textbf{print} $\tau_3$}
   \For{$n \in \tau_2$}
   \State{$\pi_{4}(n):=\pi_{4}(n)\cup (\{w\}$}
   \EndFor
   \State{$Q=\pi_{4}(x)\cap\pi_{4}(y)\cap \pi_{4}(z)\cap\pi_4(w)$}
 \For{$q \in Q$}
 \State{$\tau_4=\{x,y,z,w,q\}$}
 \State{\textbf{print} $\tau_4$}
 \EndFor
 \For{$n \in \tau_3$}
 \State{$\pi_{5}(n):=\pi_{4}(n)\cup \{q\}$}
 \EndFor
 \State{$R=\pi_{5}(x)\cap\pi_{5}(y)\cap \pi_{5}(z)\cap\pi_5(w)\cap\pi_{5}(q)$}
 \For{$r \in R$}
 \State{$\tau_5=\{x,z,y,w,q,r\}$}
 \State{\textbf{print} $\tau_5$}
 \EndFor
   \EndFor
   \EndFor
\EndFor
\end{algorithmic}
\end{algorithm}

\begin{algorithm}[H]
\caption{Compute Cliques (Recursive Evolving Dynamic Multilayered Node-Neighbourhood Approach)}
\label{alg:RDNNA}
\begin{algorithmic}[1]
    \Procedure{ComputeCliques}{$E$,$d_{\max}$}
        \State \textbf{Input:} $E$, the list of (sorted) edges, $d_{\max}$, the maximum clique dimension
        \State \textbf{Output:} All the cliques up to dimension $d_{\max}$
        \State \textbf{Initialization:}
        \State $\pi_d=\emptyset, \, \forall d \in \{2,\hdots,d_{\max}\}$
        \Statex

    \Procedure{AuxiliaryFunction}{$\sigma$}
        \State \textbf{Input:} $\sigma=(x_0,x_1,\hdots,x_{\dim})$, a clique described as an ascending sequence of nodes
        \State \textbf{Output:} A sequence of cliques $\tau$ extended from $\sigma$
        \State $d=|\sigma|$
        \State \textbf{print} $\sigma$
        \State $\hat{\sigma}:=\sigma \setminus \{x_{\dim}\}$ 
        \State
        \If{$d<d_{\max}$}

         \If{$d>2$}
         \For{ $n \in \hat{\sigma}$}
        
        \State $\pi_{d}(x):=\pi_{d}(n)\cup\{x_{\dim}\}$
        \EndFor
        \Else
         \State{$\sigma=\{x,y\}$}
         \State $\pi_{2}(x):=\pi_{2}(x)\cup\{y\}$
        \State $\pi_{2}(y):=\pi_{2}(y)\cup\{x\}$
        \EndIf
        
        \State{$Z=\cap_{n \in \sigma} \pi_d(n)$}
        \For{$z \in Z$}
        \State $\tau=\sigma \cup \{z\}$
        \State \textbf{compute} \textproc{AuxiliaryFunction}($\tau$)
        
        \EndFor
        \EndIf
    \EndProcedure

        \State
        \For{$e \in E$}
        \State \textbf{compute} \textproc{AuxiliaryFuncion}(e)
        \EndFor
        \State
    \EndProcedure
    
\end{algorithmic}
\end{algorithm}

\subsection{\md{Classic Algorithms}}
\label{sec:classic_algorithms}
The K3 algorithm \md{was developed in~\cite{chiba1985arboricity}; the forward method in~\cite{schank2005finding}}.
\begin{algorithm}[H]
\caption{Find Triangles (Edge-Iterator)}
\begin{algorithmic}[1]
\Require $E$, the List of edges; $\pi$, the node neighbourhood
\For{$e=(x,y) \in E$ }
    
   \State{$Z=\pi(x)\cap\pi(y)$}
   \For{$z \in Z$}
   \State{\textbf{print} $(x,y,v)$}
   \EndFor
\EndFor
\end{algorithmic}

\end{algorithm}

\begin{algorithm}[H]
\caption{Find Triangles (K3 Method)}
\begin{algorithmic}[1]
\Require $G=(V,E)$, an undirected simple graph
\Ensure Sorted sequence of nodes ($V=(v_1,v_2,\hdots,v_n)$) in descending order of degree
\State Create the set of marked nodes, $M=\emptyset$
\For{$i \in (1,2,\hdots,n-2)$}
\State $M:=M\cup \pi(v_i)$
\For{$u \in M$}
\For{$w \in \pi(u)$}
\If{$w \in M$}{ \textbf{print} $\{v_i,u,w\}$}
\State $M:=M\setminus \{u\}$
\EndIf
\EndFor
\EndFor
\State $V:=V\setminus\{v_i\}$
\EndFor
\end{algorithmic}

\end{algorithm}

\begin{algorithm}[H]
\caption{Find Triangles (forward method)}
\begin{algorithmic}[1]
\Require $V$, the list of nodes (ordered by node degree in decreasing order); $\pi$ the adjacency of nodes

   \State{$A(x):=\emptyset, \, \forall v \in V$}
    \For{$y \in V$}
    \For{$z \in \pi(y)$}
    \If{$y<z$}{\For{$x \in A(y)\cap A(z)$}
    \State{\textbf{print} $\{x,y,z\}$}
    \EndFor}
    \State{$A(z):=A(z)\cup\{y\}$}
    \EndIf
    \EndFor
    \EndFor

\end{algorithmic}

\end{algorithm}

\section{Tables}\label{sec:tables}
In this section, we provide the tables to elucidate the iteration process or our algorithms. We used the Figure \ref{fig:iteration_process} to provide the examples. We also provide the cliques counting from the benchmark performed in Section \ref{sec:benchmark}.
%
%
%
%
%
\begin{table}[htbp]
\hspace*{0.25cm}
\begin{tabular}{|>{\centering}p{20mm}|p{110mm}|p{35mm}|}
\hline
\textbf{Iteration} & \textbf{Clique Found} & \textbf{Boundary of Clique} \\
\hline
0 & $\emptyset$ &  \ \ \ \  \ \text{-}\\
\hline
1 & $\{1, 2\}$ & $\{\{1\}, \{2\}\}$ \\
\hline
2 & $\{2, 3\}$ & $\{\{2\}, \{3\}\}$ \\
\hline
3 & $\{1, 5\}$ & $\{\{1\}, \{5\}\}$ \\
\hline
4 & $\{2, 5\}$ & $\{\{2\}, \{5\}\}$ \\
\hline
5 & $\{1, 2, 5\}$ & $\{\{2, 5\}, \{1, 2\}, \{1, 5\}\}$ \\
\hline
6 & $\{3, 5\}$ & $\{\{3\}, \{5\}\}$ \\
\hline
7 & $\{2, 3, 5\}$ & $\{\{3, 5\}, \{2, 3\}, \{2, 5\}\}$ \\
\hline
8 & $\{1, 3\}$ & $\{\{1\}, \{3\}\}$ \\
\hline
9 & $\{1, 2, 3\}$ & $\{\{1, 3\}, \{1, 2\}, \{2, 3\}\}$ \\
\hline
10 & $\{1, 3, 5\}$ & $\{\{1, 3\}, \{1, 5\}, \{3, 5\}\}$ \\
\hline
11 & $\{1, 2, 3\}$ & $\{\{1, 3\}, \{1, 2\}, \{2, 3\}\}$ \\
\hline
12 & $\{3, 4\}$ & $\{\{3\}, \{4\}\}$ \\
\hline
13 & $\{4, 5\}$ & $\{\{5\}, \{4\}\}$ \\
\hline
14 & $\{3, 4, 5\}$ & $\{\{4, 5\}, \{3, 5\}, \{3, 4\}\}$ \\
\hline
\end{tabular}
\caption{The evolution of the cliques finding and the description of the boundary of each clique found in the iteration process of the Algorithm \ref{alg:RDBNA} on Figure \ref{fig:iteration_process}.}
\label{tab:cliques_boundaries_example}
\end{table}
\begin{table}[h]
    \hspace*{0.25cm}
    \begin{tabular}{|>{\centering}p{20mm}|p{150mm}|}
        \hline
        \textbf{Iteration} & \textbf{Updated Boundary Set} \\
        \hline
        0 & N/A \\
        \hline
        1 & $\{1\}: \{2\}, \{2\}: \{1\}$ \\
        \hline
        2 & $\{1\}: \{2\}, \{2\}: \{1, 3\}, \{3\}: \{2\}$ \\
        \hline
        3 & $\{1\}: \{2, 5\}, \{2\}: \{1, 3\}, \{3\}: \{2\}, \{5\}: \{1\}$ \\
        \hline
        4 & $\{1\}: \{2, 5\}, \{2\}: \{1, 3, 5\}, \{3\}: \{2\}, \{5\}: \{1, 2\}$ \\
        \hline
        5 & $\{1\}: \{2, 5\}, \{2\}: \{1, 3, 5\}, \{3\}: \{2\}, \{5\}: \{1, 2\}, \{2, 5\}: \{1\}, \{1, 2\}: \{5\}, \{1, 5\}: \{2\}$ \\
        \hline
        6 & $\{1\}: \{2, 5\}, \{2\}: \{1, 3, 5\}, \{3\}: \{2, 5\}, \{5\}: \{1, 2, 3\}, \{2, 5\}: \{1\}, \{1, 2\}: \{5\}, \{1, 5\}: \{2\}$ \\
        \hline
        7 & $\{1\}: \{2, 5\}, \{2\}: \{1, 3, 5\}, \{3\}: \{2, 5\}, \{5\}: \{1, 2, 3\}, \{2, 5\}: \{1, 3\}, \{1, 2\}: \{5\}, \{1, 5\}: \{2\}, \{3, 5\}: \{2\}, \{2, 3\}: \{5\}$ \\
        \hline
        8 & $\{1\}: \{2, 3, 5\}, \{2\}: \{1, 3, 5\}, \{3\}: \{1, 2, 5\}, \{5\}: \{1, 2, 3\}, \{2, 5\}: \{1, 3\}, \{1, 2\}: \{5\}, \{1, 5\}: \{2\}, \{3, 5\}: \{2\}, \{2, 3\}: \{1, 5\}$ \\
        \hline
        9 & $\{1\}: \{2, 3, 5\}, \{2\}: \{1, 3, 5\}, \{3\}: \{1, 2, 5\}, \{5\}: \{1, 2, 3\}, \{2, 5\}: \{1, 3\}, \{1, 2\}: \{3, 5\}, \{1, 5\}: \{2\}, \{3, 5\}: \{2\}, \{2, 3\}: \{1, 5\}, \{1, 3\}: \{2\}$ \\
        \hline
        10 & $\{1\}: \{2, 3, 5\}, \{2\}: \{1, 3, 5\}, \{3\}: \{1, 2, 5\}, \{5\}: \{1, 2, 3\}, \{2, 5\}: \{1, 3\}, \{1, 2\}: \{3, 5\}, \{1, 5\}: \{2, 3\}, \{3, 5\}: \{1, 2\}, \{2, 3\}: \{1, 5\}, \{1, 3\}: \{2, 5\}$ \\
        \hline
        11 & $\{1\}: \{2, 3, 5\}, \{2\}: \{1, 3, 5\}, \{3\}: \{1, 2, 5\}, \{5\}: \{1, 2, 3\}, \{2, 5\}: \{1, 3\}, \{1, 2\}: \{3, 5\}, \{1, 5\}: \{2, 3\}, \{3, 5\}: \{1, 2\}, \{2, 3\}: \{1, 5\}, \{1, 3\}: \{2, 5\}, \{1, 2, 3, 5\}$ \\
        \hline
        12 & $\{1\}: \{2, 3, 5\}, \{2\}: \{1, 3, 5\}, \{3\}: \{1, 2, 4, 5\}, \{5\}: \{1, 2, 3, 4\}, \{2, 5\}: \{1, 3\}, \{1, 2\}: \{3, 5\}, \{1, 5\}: \{2, 3\}, \{3, 5\}: \{1, 2\}, \{2, 3\}: \{1, 5\}, \{4\}: \{3\}$ \\
        \hline
        13 & $\{1\}: \{2, 3, 5\}, \{2\}: \{1, 3, 5\}, \{3\}: \{1, 2, 4, 5\}, \{5\}: \{1, 2, 3, 4\}, \{2, 5\}: \{1, 3\}, \{1, 2\}: \{3, 5\}, \{1, 5\}: \{2, 3\}, \{3, 5\}: \{1, 2\}, \{2, 3\}: \{1, 5\}, \{4\}: \{3, 5\}, \{4, 5\}: \{3\}$ \\
        \hline
        14 & $\{1\}: \{2, 3, 5\}, \{2\}: \{1, 3, 5\}, \{3\}: \{1, 2, 4, 5\}, \{5\}: \{1, 2, 3, 4\}, \{2, 5\}: \{1, 3\}, \{1, 2\}: \{3, 5\}, \{1, 5\}: \{2, 3\}, \{3, 5\}: \{1, 2\}, \{2, 3\}: \{1, 5\}, \{4\}: \{3, 5\}, \{4, 5\}: \{3\}, \{3, 4\}: \{5\}$ \\
        \hline
    \end{tabular}
    \caption{The evolution of the boundary neighbourhood of the graph in Figure \ref{fig:iteration_process}  by using the Algorithm \ref{alg:RDNNA}.}
     \label{tab:boundary-set}
\end{table}
\begin{table}[ht]
\centering
\begin{tabular}{|>{\centering}p{20mm}|p{120mm}|p{25mm}|}
\hline
\textbf{Iteration} & \textbf{Current Neighborhood Updates} & \textbf{Clique Found} \\ \hline
0  & $\pi_2$: \{1: \{\}, 2: \{\}, 3: \{\}, 5: \{\}, 4: \{\}\}, $\pi_3$: \{1: \{\}, 2: \{\}, 3: \{\}, 5: \{\}, 4: \{\}\} & - \\ \hline
1  & $\pi_2(1):= \pi_2(1) \cup \{2\} = \{2\}, \pi_2(2):= \pi_2(2) \cup \{1\} = \{1\}$ & \{1, 2\} \\ \hline
2  & $\pi_2(1):= \pi_2(1) \cup \{5\} = \{2, 5\}, \pi_2(5):= \pi_2(5) \cup \{1\} = \{1\}$ & \{1, 5\} \\ \hline
3  & $\pi_2(1):= \pi_2(1) \cup \{3\} = \{2, 3, 5\}, \pi_2(3):= \pi_2(3) \cup \{1\} = \{1\}$ & \{1, 3\} \\ \hline
4  & $\pi_2(2):= \pi_2(2) \cup \{3\} = \{1, 3\}, \pi_2(3):= \pi_2(3) \cup \{2\} = \{1, 2\}$ & \{2, 3\} \\ \hline
5  & $\pi_3(2):= \pi_3(2) \cup (\pi_3(2) - \{2, 3\}) = \{1\}, \pi_3(3):= \pi_3(3) \cup (\pi_3(3) - \{2, 3\}) = \{1\}$ & \{1, 2, 3\} \\ \hline
6  & $\pi_2(2):= \pi_2(2) \cup \{5\} = \{1, 3, 5\}, \pi_2(5):= \pi_2(5) \cup \{2\} = \{1, 2\}$ & \{2, 5\} \\ \hline
7  & $\pi_3(2):= \pi_3(2) \cup (\pi_3(2) - \{2, 5\}) = \{1\}, \pi_3(5):= \pi_3(5) \cup (\pi_3(5) - \{2, 5\}) = \{1\}$ & \{1, 2, 5\} \\ \hline
8  & $\pi_2(3):= \pi_2(3) \cup \{5\} = \{1, 2, 5\}, \pi_2(5):= \pi_2(5) \cup \{3\} = \{1, 2, 3\}$ & \{3, 5\} \\ \hline
9  & $\pi_3(3):= \pi_3(3) \cup (\pi_3(3) - \{3, 5\}) = \{1\}, \pi_3(5):= \pi_3(5) \cup (\pi_3(5) - \{3, 5\}) = \{1\}$ & \{1, 3, 5\} \\ \hline
10 & $\pi_3(3):= \pi_3(3) \cup (\pi_3(3) - \{3, 5\}) = \{1, 2\}, \pi_3(5):= \pi_3(5) \cup (\pi_3(5) - \{3, 5\}) = \{1, 2\}$ & \{2, 3, 5\} \\ \hline
11 & No neighborhood updates. & \{1, 2, 3, 5\} \\ \hline
12 & $\pi_2(3):= \pi_2(3) \cup \{4\} = \{1, 2, 4, 5\}, \pi_2(4):= \pi_2(4) \cup \{3\} = \{3\}$ & \{3, 4\} \\ \hline
13 & $\pi_2(5):= \pi_2(5) \cup \{4\} = \{1, 2, 3, 4\}, \pi_2(4):= \pi_2(4) \cup \{5\} = \{3, 5\}$ & \{4, 5\} \\ \hline
14 & $\pi_3(5):= \pi_3(5) \cup (\pi_3(5) - \{4, 5\}) = \{1, 2, 3\}, \pi_3(4):= \pi_3(4) \cup (\pi_3(4) - \{4, 5\}) = \{3\}$ & \{3, 4, 5\} \\ \hline
\end{tabular}
\caption{The evolution of the cliques search in the example of Figure \ref{fig:iteration_process} by using the Evolving Dynamic Multilayered Node-Neighbourhood Approach (Algorithm \ref{alg:RDNNA}).}
\label{tab:node-neigh}
\end{table}
%
%
\begin{table}[h!]
\centering
\begin{tabular}{|l|l|r|r|r|l|}
\hline
 & \textbf{Method} & $d_{\max}$ & \textbf{Number of nodes} & \textbf{Number of edges} & \textbf{Number of cliques} \\
\hline
0 & \text{NetworkX} & 2 & 50515 & 819090 & 3143305 \\
\hline
1 & \text{Algorithm 9} & 2 & 50515 & 819090 & 3143305 \\
\hline
2 & \text{Algorithm 5} & 2 & 50515 & 819090 & 3143305 \\
\hline
3 & \text{NetworkX} & 3 & 50515 & 819090 & 7613489 \\
\hline
4 & \text{Algorithm 9} & 3 & 50515 & 819090 & 7613489 \\
\hline
5 & \text{Algorithm 5} & 3 & 50515 & 819090 & 7613489 \\
\hline
6 & \text{NetworkX} & 4 & 50515 & 819090 & 14757799 \\
\hline
7 & \text{Algorithm 9} & 4 & 50515 & 819090 & 14757799 \\
\hline
8 & \text{Algorithm 5} & 4 & 50515 & 819090 & 14757799 \\
\hline
9 & \text{NetworkX} & 5 & 50515 & 819090 & 24121489 \\
\hline
10 & \text{Algorithm 9} & 5 & 50515 & 819090 & 24121489 \\
\hline
11 & \text{Algorithm 5} & 5 & 50515 & 819090 & 24121489 \\
\hline
12 & \text{Algorithm 5} & 6 & 50515 & 819090 & 34313971 \\
\hline
13 & \text{Algorithm 9} & 6 & 50515 & 819090 & 34313971 \\
\hline
14 & \text{NetworkX} & 6 & 50515 & 819090 & 34313971 \\
\hline
15 & \text{NetworkX} & 7 & 50515 & 819090 & \text{Not computed} \\
\hline
16 & \text{Algorithm 9} & 7 & 50515 & 819090 & 43616877 \\
\hline
17 & \text{Algorithm 5} & 7 & 50515 & 819090 & 43616877 \\
\hline
18 & \text{NetworkX} & 8 & 50515 & 819090 & \text{Not computed} \\
\hline
19 & \text{Algorithm 9} & 8 & 50515 & 819090 & 50728110 \\
\hline
20 & \text{Algorithm 5} & 8 & 50515 & 819090 & \text{Not computed} \\
\hline
21 & \text{NetworkX} & 9 & 50515 & 819090 & \text{Not computed} \\
\hline
22 & \text{Algorithm 9} & 9 & 50515 & 819090 & 55232908 \\
\hline
23 & \text{Algorithm 5} & 9 & 50515 & 819090 & \text{Not computed} \\
\hline
24 & \text{NetworkX} & 10 & 50515 & 819090 & \text{Not computed} \\
\hline
25 & \text{Algorithm 9} & 10 & 50515 & 819090 & 57568864 \\
\hline
26 & \text{Algorithm 5} & 10 & 50515 & 819090 & \text{Not computed} \\
\hline
\end{tabular}
\caption{Comparison of different methods for computing cliques from the file \bf{ artist\_edges.csv} .}
\label{table:cliques_count_artist_edges}
\end{table}
\begin{table}[h!]
\centering
\begin{tabular}{|l|l|r|r|r|l|}
\hline
 & \textbf{Method} & $d_{\max}$ & \textbf{Number of nodes} & \textbf{Number of edges} & \textbf{Number of cliques} \\
\hline
0 & \text{NetworkX} & 2 & 3892 & 17239 & 108221 \\
\hline
1 & \text{Algorithm 9} & 2 & 3892 & 17239 & 108221 \\
\hline
2 & \text{Algorithm 5} & 2 & 3892 & 17239 & 108221 \\
\hline
3 & \text{NetworkX} & 3 & 3892 & 17239 & 904252 \\
\hline
4 & \text{Algorithm 9} & 3 & 3892 & 17239 & 904252 \\
\hline
5 & \text{Algorithm 5} & 3 & 3892 & 17239 & 904252 \\
\hline
6 & \text{NetworkX} & 4 & 3892 & 17239 & \text{Not computed} \\
\hline
7 & \text{Algorithm 9} & 4 & 3892 & 17239 & 8465416 \\
\hline
8 & \text{Algorithm 5} & 4 & 3892 & 17239 & 8465416 \\
\hline
9 & \text{Algorithm 5} & 5 & 3892 & 17239 & 73114297 \\
\hline
10 & \text{Algorithm 9} & 5 & 3892 & 17239 & 73114297 \\
\hline
11 & \text{NetworkX} & 5 & 3892 & 17239 & \text{Not computed} \\
\hline
12 & \text{NetworkX} & 6 & 3892 & 17239 & \text{Not computed} \\
\hline
13 & \text{Algorithm 9} & 6 & 3892 & 17239 & 555455386 \\
\hline
14 & \text{Algorithm 5} & 6 & 3892 & 17239 & 555455386 \\
\hline
15 & \text{NetworkX} & 7 & 3892 & 17239 & \text{Not computed} \\
\hline
16 & \text{Algorithm 9} & 7 & 3892 & 17239 & 3686393585 \\
\hline
17 & \text{Algorithm 5} & 7 & 3892 & 17239 & \text{Not computed} \\
\hline
18 & \text{NetworkX} & 8 & 3892 & 17239 & \text{Not computed} \\
\hline
19 & \text{Algorithm 9} & 8 & 3892 & 17239 & 21466537882 \\
\hline
20 & \text{Algorithm 5} & 8 & 3892 & 17239 & \text{Not computed} \\
\hline
\end{tabular}
\caption{Comparison of different methods for computing cliques from the file \bf{ tvshow\_edges.csv} .}
\label{table:cliques_count_tvshow_edges}
\end{table}
\begin{table}[h!]
\centering
\begin{tabular}{|l|l|r|r|r|l|}
\hline
 & \textbf{Method} & $d_{\max}$ & \textbf{Number of nodes} & \textbf{Number of edges} & \textbf{Number of cliques} \\
\hline
0 & \text{Algorithm 5} & 2 & 14113 & 52126 & 122244 \\
\hline
1 & \text{NetworkX} & 2 & 14113 & 52126 & 122244 \\
\hline
2 & \text{Algorithm 9} & 2 & 14113 & 52126 & 122244 \\
\hline
3 & \text{Algorithm 5} & 3 & 14113 & 52126 & 222554 \\
\hline
4 & \text{NetworkX} & 3 & 14113 & 52126 & 222554 \\
\hline
5 & \text{Algorithm 9} & 3 & 14113 & 52126 & 222554 \\
\hline
6 & \text{Algorithm 5} & 4 & 14113 & 52126 & 430383 \\
\hline
7 & \text{NetworkX} & 4 & 14113 & 52126 & 430383 \\
\hline
8 & \text{Algorithm 9} & 4 & 14113 & 52126 & 430383 \\
\hline
9 & \text{Algorithm 5} & 5 & 14113 & 52126 & 834921 \\
\hline
10 & \text{NetworkX} & 5 & 14113 & 52126 & 834921 \\
\hline
11 & \text{Algorithm 9} & 5 & 14113 & 52126 & 834921 \\
\hline
12 & \text{Algorithm 9} & 6 & 14113 & 52126 & 1536280 \\
\hline
13 & \text{NetworkX} & 6 & 14113 & 52126 & 1536280 \\
\hline
14 & \text{Algorithm 5} & 6 & 14113 & 52126 & 1536280 \\
\hline
15 & \text{NetworkX} & 6 & 14113 & 52126 & 1536280 \\
\hline
16 & \text{Algorithm 5} & 7 & 14113 & 52126 & 2600400 \\
\hline
17 & \text{Algorithm 5} & 7 & 14113 & 52126 & 2600400 \\
\hline
18 & \text{NetworkX} & 7 & 14113 & 52126 & \text{Not computed} \\
\hline
19 & \text{Algorithm 9} & 7 & 14113 & 52126 & 2600400 \\
\hline
20 & \text{Algorithm 5} & 8 & 14113 & 52126 & \text{Not computed} \\
\hline
21 & \text{NetworkX} & 8 & 14113 & 52126 & \text{Not computed} \\
\hline
22 & \text{Algorithm 9} & 8 & 14113 & 52126 & 3994556 \\
\hline
23 & \text{Algorithm 5} & 9 & 14113 & 52126 & \text{Not computed} \\
\hline
24 & \text{NetworkX} & 9 & 14113 & 52126 & \text{Not computed} \\
\hline
25 & \text{Algorithm 9} & 9 & 14113 & 52126 & 5551182 \\
\hline
26 & \text{NetworkX} & 10 & 14113 & 52126 & \text{Not computed} \\
\hline
27 & \text{Algorithm 5} & 10 & 14113 & 52126 & \text{Not computed} \\
\hline
28 & \text{Algorithm 9} & 10 & 14113 & 52126 & 7016890 \\
\hline
\end{tabular}
\caption{Comparison of different methods for computing cliques from the file \bf{ company\_edges.csv} .}
\label{table:cliques_count_company_edges}
\end{table}
\begin{table}[h!]
\centering
\begin{tabular}{|l|l|r|r|r|l|}
\hline
 & \textbf{Method} & $d_{\max}$ & \textbf{Number of nodes} & \textbf{Number of edges} & \textbf{Number of cliques} \\
\hline
0 & \text{NetworkX} & 2 & 11565 & 67038 & 263441 \\
\hline
1 & \text{Algorithm 9} & 2 & 11565 & 67038 & 263441 \\
\hline
2 & \text{Algorithm 5} & 2 & 11565 & 67038 & 263441 \\
\hline
3 & \text{NetworkX} & 3 & 11565 & 67038 & 837273 \\
\hline
4 & \text{Algorithm 9} & 3 & 11565 & 67038 & 837273 \\
\hline
5 & \text{Algorithm 5} & 3 & 11565 & 67038 & 837273 \\
\hline
6 & \text{NetworkX} & 4 & 11565 & 67038 & 2228890 \\
\hline
7 & \text{Algorithm 9} & 4 & 11565 & 67038 & 2228890 \\
\hline
8 & \text{Algorithm 5} & 4 & 11565 & 67038 & 2228890 \\
\hline
9 & \text{NetworkX} & 5 & 11565 & 67038 & 4902092 \\
\hline
10 & \text{Algorithm 9} & 5 & 11565 & 67038 & 4902092 \\
\hline
11 & \text{Algorithm 5} & 5 & 11565 & 67038 & 4902092 \\
\hline
12 & \text{Algorithm 5} & 6 & 11565 & 67038 & 9336183 \\
\hline
13 & \text{Algorithm 9} & 6 & 11565 & 67038 & 9336183 \\
\hline
14 & \text{NetworkX} & 6 & 11565 & 67038 & 9336183 \\
\hline
15 & \text{NetworkX} & 7 & 11565 & 67038 & \text{Not computed} \\
\hline
16 & \text{Algorithm 9} & 7 & 11565 & 67038 & 16045175 \\
\hline
17 & \text{Algorithm 5} & 7 & 11565 & 67038 & 16045175 \\
\hline
18 & \text{NetworkX} & 8 & 11565 & 67038 & \text{Not computed} \\
\hline
19 & \text{Algorithm 9} & 8 & 11565 & 67038 & 25276956 \\
\hline
20 & \text{Algorithm 5} & 8 & 11565 & 67038 & \text{Not computed} \\
\hline
21 & \text{NetworkX} & 9 & 11565 & 67038 & \text{Not computed} \\
\hline
22 & \text{Algorithm 9} & 9 & 11565 & 67038 & 36454343 \\
\hline
23 & \text{Algorithm 5} & 9 & 11565 & 67038 & \text{Not computed} \\
\hline
24 & \text{NetworkX} & 10 & 11565 & 67038 & \text{Not computed} \\
\hline
25 & \text{Algorithm 9} & 10 & 11565 & 67038 & 48016033 \\
\hline
26 & \text{Algorithm 5} & 10 & 11565 & 67038 & \text{Not computed} \\
\hline
\end{tabular}
\caption{Comparison of different methods for computing cliques from the file \bf{ public\_figure\_edges.csv} .}
\label{table:cliques_count_public_figure_edges}
\end{table}
\begin{table}[h!]
\centering
\begin{tabular}{|l|l|r|r|r|l|}
\hline
 & \textbf{Method} & $d_{\max}$ & \textbf{Number of nodes} & \textbf{Number of edges} & \textbf{Number of cliques} \\
\hline
0 & \text{Algorithm 9} & 2 & 7057 & 89429 & 620340 \\
\hline
1 & \text{Algorithm 5} & 2 & 7057 & 89429 & 620340 \\
\hline
2 & \text{NetworkX} & 2 & 7057 & 89429 & 620340 \\
\hline
3 & \text{Algorithm 9} & 3 & 7057 & 89429 & 2885484 \\
\hline
4 & \text{Algorithm 5} & 3 & 7057 & 89429 & 2885484 \\
\hline
5 & \text{NetworkX} & 3 & 7057 & 89429 & 2885484 \\
\hline
6 & \text{Algorithm 9} & 4 & 7057 & 89429 & 10357710 \\
\hline
7 & \text{Algorithm 5} & 4 & 7057 & 89429 & 10357710 \\
\hline
8 & \text{NetworkX} & 4 & 7057 & 89429 & 10357710 \\
\hline
9 & \text{Algorithm 9} & 5 & 7057 & 89429 & 30044661 \\
\hline
10 & \text{Algorithm 5} & 5 & 7057 & 89429 & 30044661 \\
\hline
11 & \text{NetworkX} & 5 & 7057 & 89429 & 30044661 \\
\hline
12 & \text{NetworkX} & 6 & 7057 & 89429 & \text{Not computed} \\
\hline
13 & \text{Algorithm 5} & 6 & 7057 & 89429 & 73763797 \\
\hline
14 & \text{Algorithm 9} & 6 & 7057 & 89429 & 73763797 \\
\hline
15 & \text{Algorithm 9} & 7 & 7057 & 89429 & 159818173 \\
\hline
16 & \text{Algorithm 5} & 7 & 7057 & 89429 & \text{Not computed} \\
\hline
17 & \text{NetworkX} & 7 & 7057 & 89429 & \text{Not computed} \\
\hline
18 & \text{Algorithm 9} & 8 & 7057 & 89429 & 314967335 \\
\hline
19 & \text{Algorithm 5} & 8 & 7057 & 89429 & \text{Not computed} \\
\hline
20 & \text{NetworkX} & 8 & 7057 & 89429 & \text{Not computed} \\
\hline
21 & \text{Algorithm 9} & 9 & 7057 & 89429 & 573243192 \\
\hline
22 & \text{Algorithm 5} & 9 & 7057 & 89429 & \text{Not computed} \\
\hline
23 & \text{NetworkX} & 9 & 7057 & 89429 & \text{Not computed} \\
\hline
24 & \text{Algorithm 9} & 10 & 7057 & 89429 & 966494969 \\
\hline
25 & \text{Algorithm 5} & 10 & 7057 & 89429 & \text{Not computed} \\
\hline
26 & \text{NetworkX} & 10 & 7057 & 89429 & \text{Not computed} \\
\hline
\end{tabular}
\caption{Comparison of different methods for computing cliques from the file \bf{ government\_edges.csv} .}
\label{table:cliques_count_government_edges}
\end{table}
\begin{table}[h!]
\centering
\begin{tabular}{|l|l|r|r|r|l|}
\hline
 & \textbf{Method} & $d_{\max}$ & \textbf{Number of nodes} & \textbf{Number of edges} & \textbf{Number of cliques} \\
\hline
0 & \text{Algorithm 5} & 2 & 5908 & 41706 & 222246 \\
\hline
1 & \text{NetworkX} & 2 & 5908 & 41706 & 222246 \\
\hline
2 & \text{Algorithm 9} & 2 & 5908 & 41706 & 222246 \\
\hline
3 & \text{Algorithm 5} & 3 & 5908 & 41706 & 877240 \\
\hline
4 & \text{NetworkX} & 3 & 5908 & 41706 & 877240 \\
\hline
5 & \text{Algorithm 9} & 3 & 5908 & 41706 & 877240 \\
\hline
6 & \text{NetworkX} & 4 & 5908 & 41706 & 2879490 \\
\hline
7 & \text{Algorithm 9} & 4 & 5908 & 41706 & 2879490 \\
\hline
8 & \text{Algorithm 5} & 4 & 5908 & 41706 & 2879490 \\
\hline
9 & \text{Algorithm 5} & 5 & 5908 & 41706 & 7800497 \\
\hline
10 & \text{NetworkX} & 5 & 5908 & 41706 & 7800497 \\
\hline
11 & \text{Algorithm 9} & 5 & 5908 & 41706 & 7800497 \\
\hline
12 & \text{Algorithm 5} & 6 & 5908 & 41706 & 17627458 \\
\hline
13 & \text{NetworkX} & 6 & 5908 & 41706 & 17627458 \\
\hline
14 & \text{Algorithm 9} & 6 & 5908 & 41706 & 17627458 \\
\hline
15 & \text{NetworkX} & 7 & 5908 & 41706 & \text{Not computed} \\
\hline
16 & \text{Algorithm 5} & 7 & 5908 & 41706 & 33781950 \\
\hline
17 & \text{Algorithm 9} & 7 & 5908 & 41706 & 33781950 \\
\hline
\end{tabular}
\caption{Comparison of different methods for computing cliques from the file \bf{ politician\_edges.csv} .}
\label{table:cliques_count_politician_edges}
\end{table}
\begin{table}[h!]
\centering
\begin{tabular}{|l|l|r|r|r|l|}
\hline
 & \textbf{Method} & $d_{\max}$ & \textbf{Number of nodes} & \textbf{Number of edges} & \textbf{Number of cliques} \\
\hline
0  & \text{Algorithm 5} & 2  & 13866 & 86811 & 240700 \\
\hline
1  & \text{Algorithm 9} & 2  & 13866 & 86811 & 240700 \\
\hline
2  & \text{NetworkX}    & 2  & 13866 & 86811 & 240700 \\
\hline
3  & \text{Algorithm 5} & 3  & 13866 & 86811 & 500893 \\
\hline
4  & \text{Algorithm 9} & 3  & 13866 & 86811 & 500893 \\
\hline
5  & \text{NetworkX}    & 3  & 13866 & 86811 & 500893 \\
\hline
6  & \text{Algorithm 5} & 4  & 13866 & 86811 & 1364416 \\
\hline
7  & \text{Algorithm 9} & 4  & 13866 & 86811 & 1364416 \\
\hline
8  & \text{NetworkX}    & 4  & 13866 & 86811 & 1364416 \\
\hline
9  & \text{Algorithm 5} & 5  & 13866 & 86811 & 4666654 \\
\hline
10 & \text{Algorithm 9} & 5  & 13866 & 86811 & 4666654 \\
\hline
11 & \text{NetworkX}    & 5  & 13866 & 86811 & 4666654 \\
\hline
12 & \text{NetworkX}    & 6  & 13866 & 86811 & 15922951 \\
\hline
13 & \text{Algorithm 9} & 6  & 13866 & 86811 & 15922951 \\
\hline
14 & \text{Algorithm 5} & 6  & 13866 & 86811 & 15922951 \\
\hline
15 & \text{Algorithm 5} & 7  & 13866 & 86811 & 48663702 \\
\hline
16 & \text{Algorithm 9} & 7  & 13866 & 86811 & 48663702 \\
\hline
17 & \text{NetworkX}    & 7  & 13866 & 86811 & \text{Not computed} \\
\hline
18 & \text{Algorithm 5} & 8  & 13866 & 86811 & \text{Not computed} \\
\hline
19 & \text{Algorithm 9} & 8  & 13866 & 86811 & 129868129 \\
\hline
20 & \text{NetworkX}    & 8  & 13866 & 86811 & \text{Not computed} \\
\hline
21 & \text{Algorithm 5} & 9  & 13866 & 86811 & \text{Not computed} \\
\hline
22 & \text{Algorithm 9} & 9  & 13866 & 86811 & 302555642 \\
\hline
23 & \text{NetworkX}    & 9  & 13866 & 86811 & \text{Not computed} \\
\hline
24 & \text{Algorithm 5} & 10 & 13866 & 86811 & \text{Not computed} \\
\hline
25 & \text{Algorithm 9} & 10 & 13866 & 86811 & 619197827 \\
\hline
26 & \text{NetworkX}    & 10 & 13866 & 86811 & \text{Not computed} \\
\hline
\end{tabular}
\caption{Comparison of different methods for computing cliques from the file \bf{ athletes\_edges.csv} .}
\label{table:cliques_count_athletes_edges}
\end{table}
\begin{table}[h!]
\centering
\begin{tabular}{|l|l|r|r|r|l|}
\hline
 & \textbf{Method} & $d_{\max}$ & \textbf{Number of nodes} & \textbf{Number of edges} & \textbf{Number of cliques} \\
\hline
0  & \text{NetworkX}    & 2  & 27917 & 205964 & 621325 \\
\hline
1  & \text{Algorithm 9} & 2  & 27917 & 205964 & 621325 \\
\hline
2  & \text{Algorithm 5} & 2  & 27917 & 205964 & 621325 \\
\hline
3  & \text{NetworkX}    & 3  & 27917 & 205964 & 1440079 \\
\hline
4  & \text{Algorithm 9} & 3  & 27917 & 205964 & 1440079 \\
\hline
5  & \text{Algorithm 5} & 3  & 27917 & 205964 & 1440079 \\
\hline
6  & \text{Algorithm 9} & 4  & 27917 & 205964 & 3583959 \\
\hline
7  & \text{Algorithm 5} & 4  & 27917 & 205964 & 3583959 \\
\hline
8  & \text{NetworkX}    & 4  & 27917 & 205964 & 3583959 \\
\hline
9  & \text{NetworkX}    & 5  & 27917 & 205964 & 9584421 \\
\hline
10 & \text{Algorithm 9} & 5  & 27917 & 205964 & 9584421 \\
\hline
11 & \text{Algorithm 5} & 5  & 27917 & 205964 & 9584421 \\
\hline
12 & \text{NetworkX}    & 6  & 27917 & 205964 & 25736067 \\
\hline
13 & \text{Algorithm 9} & 6  & 27917 & 205964 & 25736067 \\
\hline
14 & \text{Algorithm 5} & 6  & 27917 & 205964 & 25736067 \\
\hline
15 & \text{Algorithm 9} & 7  & 27917 & 205964 & 65328247 \\
\hline
16 & \text{NetworkX}    & 7  & 27917 & 205964 & \text{Not computed} \\
\hline
17 & \text{Algorithm 5} & 7  & 27917 & 205964 & 65328247 \\
\hline
\end{tabular}
\caption{Comparison of different methods for computing cliques from the file \bf{ new\_sites\_edges.csv} .}
\label{table:cliques_count_new_sites_edges}
\end{table}

\clearpage
\bibliographystyle{plain}
\bibliography{references}
\end{document}